\documentclass[10pt,a4paper]{amsart}
\usepackage{amsmath}
\usepackage[utf8]{inputenc}
\usepackage{amssymb,amsmath}
\usepackage{amsthm}
\usepackage{epsfig}
\usepackage{graphicx}
\usepackage{bm}
\usepackage{color}
\usepackage{bm}
\usepackage{amsrefs}
\usepackage{pdfpages}
\usepackage{lineno}
\usepackage{setspace} 
\usepackage{mathrsfs}
\input xypic
\usepackage{xcolor}
\usepackage{soul}
\usepackage{todonotes}
\setuptodonotes{inline}
\usepackage{mathtools}
\input xypic
\usepackage{vmargin}
\usepackage{ulem}
\setpapersize{A4}
\setmargins{2.5cm}       
{1.5cm}                        
{16.5cm}                      
{23.42cm}                    
{10pt}                           
{1cm}                           
{0pt}                             
{2cm}                           

\newtheorem{theorem}{Theorem}[section]
\newtheorem{definition}[theorem]{Definition}
\newtheorem{lemma}[theorem]{Lemma}
\newtheorem{corollary}[theorem]{Corollary}
\newtheorem{proposition}[theorem]{Proposition}
\newtheorem*{stheorem}{Theorem}
\theoremstyle{definition}
\newtheorem{remark}{Remark}[section]

\title{Soft random simplicial complexes } \author{Juli\'an David Candela
  Coca} 
  
\begin{document}
\vspace{0.5cm}
\large

\begin{abstract}
  \large A soft random graph $G(n,r,p)$ can be obtained from the
  random geometric graph $G(n,r)$ by keeping every edge in $G(n,r)$
  with probability $p$. This random graph is a particular case of the
  soft random graph model introduced by Penrose in \cite{MR3476631},
  in which the probability between two vertices is a function that
  depends on the distance between them. In this article, we define
  models for random simplicial complexes built over the soft random
  graph $G(n,r,p)$, which also present randomness in higher dimensions. Furthermore, we study the homology of those random simplicial complexes in different regimes of $n,r$ and $p$ by
  giving asymptotic formulas for the expectation of the Betti
  numbers in the sparser regimes and bounds in the denser regimes.
\end{abstract}

\maketitle

\section{Introduction}

In recent years, we have seen an
increased interest in models of random topological spaces, which has
its beginnings in the study of random graphs. The random
geometric graph $G(n,r)$, which is the point of departure for the present study, is a random graph over a
set of $n$ uniformly distributed points on the set $\mathbb{R}^d$, in
which an edge is included if the distance between its vertices is less
than $r$. In \cite{MR1986198}, Penrose extensively studied the
properties of this random graph in three different regimes: The
subcritical regime, where the geometric graph is mostly disconnected;
the critical regime, where a giant component first appears and
percolation occurs; and the supercritical regime, where finally the
geometric graph becomes connected.

A number of random simplicial complexes may be built over random
geometric graphs. These models, such as the random Vietoris-Rips
complex $R(n,r)$ and the random \v{C}ech complex $C(n,r)$, can be viewed
as higher dimensional generalizations of the random geometric graph
$G(n,r)$. In \cite{MR2770552}, Kahle studied the homology of these two
random geometric simplicial complexes in 4 regimes: The subcritical
regime, where homology first appears; the critical regime, where the
expected value of the Betti numbers is asymptotically bounded above by
$n$; the supercritial regime, where the expected value of the Betti
numbers of $R(n,r)$ dominated by $n$, and the connected regime, where homology vanishes
again.

In \cite{MR3509567}, Costa and Farber described a model for random simplicial complexes, built on a set on $n$ points, that depend on a probability vector  $\rho=(p_0,\dots,p_r)$ which allows the model to have randomness in every possible dimension. Costa and Farber, in their series of papers \cite{MR3509567},\cite{MR3661651},\cite{MR3604492}, studied important topological properties such as  connectivity, the fundamental group, and the behavior of the Betti numbers for this model.
More recently, in \cite{MR4252193}, the authors extended this research by proving limit theorems for topological invariants in  a dynamic version of Costa and Farber's model.

Soft random simplicial complexes $R(n,r,\rho)$ and
$C(n,r,\rho)$, which we study in this paper, are a generalization of
the random geometric complexes $R(n,r)$ and $C(n,r)$ that depend on a probability vector $\rho=(p_1,p_2,\dots)$. In this model,
each simplicial complex is obtained from the original one by keeping
each $k$-face with probability $p_k$. 
Thus, they can be seen as simplicial complexes built over the soft
random graph $G(n,r,p)$ \cite{MR3476631} or as a
$(p,q)$-perturbed random geometric graph $G_n^{q,p}$ \cites{MR4042102, MR4605137}, taking $p=1-p_1$ and
$q=0$. The principal contribution of this paper is the study of the
higher homology $H_k$, $k>0$ for $R(n,r,\rho)$ and
$C(n,r,\rho)$. For that, we give asymptotic bounds for the
expectation of the Betti numbers in each regime. We also prove that, similarly to the homology of
$R(n,r)$ and $C(n,r)$, the higher homology for the soft random
graphs turns out to be non-monotone.

The paper is organized as follows. The remainder of
this section will be used to introduce notation, state several
results from other references which we will use in the following, and state our main results. In Section
\ref{section2}, we define the soft random simplicial complexes
$R(n,r,\rho)$ and $C(n,r,\rho)$ and introduce the tools we need for
calculations. The remaining sections will be used to study the
behavior of the model in different regimes. The subcritical regime
will be studied in Section \ref{section3}. In this
regime, the homology for both $R(n,r,\rho)$ and $C(n,r,\rho)$ first
appear. In Section \ref{section4}, we address the critical regime and
prove that the expected value of the Betti numbers is asymptotically
bounded above by $n$, up to a constant factor for both soft random
graphs. Section \ref{section5} will be used to study the supercritical
regime. Here, we prove that the expected value of the Betti numbers of
$R(n,r,\rho)$ is asymptotically dominated by $n$. Finally, in Section
\ref{section6} we prove that, under the right assumptions,
$R(n,r,\rho)$ is asymptotically $k$-connected while $C(n,r,\rho)$ is
contractible a.a.s. Furthermore, homology for both soft random
simplicial complexes vanishes again.

The results contained in this article form a part of my
Ph.D. research at the Centro de Investigaci\'on en
Matem\'aticas, A.C. under the supervision of Antonio Rieser.

\subsection{Notation}
The random geometric graph  $G(n,r)$ was  extensively studied by Penrose in his book Random geometric graphs \cite{MR1986198}. 

Given $n$ independent, identically distributed random variables  $X_1,\dots,X_n$ in $\mathbb{R}^d$ having a common probability density function $f:\mathbb{R}^d\rightarrow [0,\infty)$, the point process $\mathcal{X}_n$ is the union $\mathcal{X}_n=\cup_{i=1}^n\{X_i\}$. Through this document, we assume that $f$  satisfies
$$\int_{\mathbb{R}^d}f(x)dx=1$$
and that for every $A\subset \mathbb{R}^d$ we have 
$$P(X_1\in A)=\int_{A}f(x)dx.$$
The random geometric graph $G(n,r)$, also denoted by $G(\mathcal{X}_n, r)$ is the the random  graph with vertex set $\mathcal{X}_n=\{X_1,\dots,X_n\}$ and an edge between $X_i$ and $X_j$  when $d(X_i, X_j)\leq r$. This graph was thoroughly studied by Penrose in \cite{MR1986198}.
On the other hand, we state the following coupling presented in
\cite{MR1986198}.  Given $\lambda>0$, let $N_\lambda$ be a Poisson
random variable independent of $\mathcal{X}_n=\{X_1,\dots,X_n\}$ and let
$$\mathcal{P}_\lambda= \{X_1,\dots,X_{N_\lambda}\}.$$
The following result from \cite{MR1986198} guarantees
that $\mathcal{P}_\lambda$ is in fact a Poisson process.
\begin{proposition}[Propositions 1.5 in \cite{MR1986198}]
$\mathcal{P}_\lambda$ is a Poisson process on $\mathbb{R}^d$ with intensity $\lambda f$.
\end{proposition}

We now define some notation from \cite{MR1986198} that
we will use along this document. Given a finite point set
$Y\subset \mathbb{R}^d$, the first element of $Y$ according to the
lexicographic order will be called the left most point of $Y$ and will
be denoted by $LMP(Y)$. For an open set $A\subset \mathbb{R}^d$  with Lebesgue measure of its boundary $Leb(\partial A)=0$  and a connected graph $\Gamma$ on $k$ vertices,  let $G_{n,A}(\Gamma)$ and  $G'_{n,A}(\Gamma)$ be the random variables  that count the number of induced subgraphs of $G(n,r)$ isomorphic to $\Gamma$ and  the number of induced subgraphs of $G(\mathcal{P}_n,r)$ isomorphic to $\Gamma$, respectively for which the left most point of the vertex set lies in $A$. Similarly, for $\Gamma$ as above,  let $J_{n,A}(\Gamma)$ and $J'_{n,A}(\Gamma)$ be the random variables  that count the number of components of $G(n,r)$ isomorphic to $\Gamma$ and  the number components of $G(\mathcal{P}_n,r)$ isomorphic to $\Gamma$, respectively, for which the left most point of the vertex set lies in $A$.

On the other hand, given a
  connected graph $\Gamma$  of order $k$ and
sets $Y$ and $A$ as above, define the 
characteristic functions $h_\Gamma$ as:
\begin{equation}
h_{\Gamma}(Y)=1_{\{G(Y,1)\cong \Gamma\}},~~~h_{\Gamma,n,A}(Y)=1_{\{G(Y,r)\cong \Gamma\}\cap \{LMP(Y)\in A\}}
\end{equation}
and $\mu_{\Gamma,A}$ as:
\begin{equation}\label{mu}
\mu_{\Gamma, A}=\frac{1}{k!}\int_A f(x)^k~dx \int_{(\mathbb{R}^d)^{k-1}}h_{\Gamma}(\{0,x_1,\dots,x_{k-1}\})~dx_1dx_2\dots dx_{k-1}
\end{equation}
notice that $\mu_{\Gamma, A}$ is finite since it is bounded by a constant times $||f||Leb(B(0,k))$. In case $A=\mathbb{R}^d$, the  subscript $A$ in all  the notation stated above will be omitted. Also, when $\Gamma=K_m$ is the complete graph on $m$ vertices, we write $\mu_{m, A}$ for $\mu_{\Gamma, A}=\mu_{K_m, A}$.

We will now state some important results from \cite{MR1986198}:   

\begin{definition}
Suppose that $\Gamma$ is a graph of order $k\geq 2$. We say that $\Gamma$ is feasible iff $P(G(\mathcal{X}_k,r)\cong \Gamma)>0$ for some $r>0$.
\end{definition}

\begin{proposition}[Propositions 3.1 in \cite{MR1986198}]\label{Penrose-1}
Suppose that $\Gamma$ is a feasible connected graph of order $k\geq 2$ and $A\subset \mathbb{R}^d$ is open with $Leb(\partial A)=0$. Also, suppose that $r=r(n)\rightarrow 0$. Then
$$\lim_{n\rightarrow \infty}\frac{G_{n,A}(\Gamma)}{n^k r^{d(k-1)}}=\lim_{n\rightarrow \infty}\frac{G_{n,A}'(\Gamma)}{n^k r^{d(k-1)}}=\mu_{\Gamma,A},$$
where $\mu_{\Gamma,A}$ is the constant defined in (\ref{mu}),  which depends on $\Gamma$ and $k$.
\end{proposition}

\begin{proposition}[Propositions 3.2 in \cite{MR1986198}]\label{Penrose-2}
Let $\Gamma$ be a feasible connected graph of order $k\geq 2$ and let $J_{n,A}(\Gamma)$ be the number of components of $G(n,r)$ isomorphic to $\Gamma$ with LMP of the set of vertices in $A$. Also, suppose that $nr^d\rightarrow 0$. Then
$$\lim_{n\rightarrow \infty}\frac{J_{n,A}(\Gamma)}{n^k r^{d(k-1)}}=\lim_{n\rightarrow \infty}\frac{J'_{n,A}(\Gamma)}{n^k r^{d(k-1)}}=\mu_{\Gamma,A}$$
where $\mu_\Gamma$ a constant that depends on $\Gamma$ and $k$.
\end{proposition}
\begin{proposition}[Propositions 3.3 in \cite{MR1986198}]\label{Penrose-3}
Suppose that $\Gamma$ is a feasible connected graph of order $k\geq 2$. Also, suppose that $nr^d\rightarrow \lambda\in(0,\infty)$. Then
$$\lim_{n\rightarrow \infty}\frac{J_{n,A}(\Gamma)}{n}=\lim_{n\rightarrow \infty}\frac{J'_{n,A}(\Gamma)}{n}=c_{\Gamma, A}$$
where $c_{\Gamma, A}$ a constant that depends on $\Gamma$ and $A$.
\end{proposition}
Furthermore, it is noticed by Kahle in \cite{MR2770552}, that Proposition~\ref{Penrose-3} can be extended to work for $G_{n,A}(\Gamma )$ as follows:
\begin{proposition}\label{Kahle-1}
Suppose that $\Gamma$ is a feasible connected graph of order $k\geq 2$. Also, suppose that $nr^d\rightarrow \lambda\in(0,\infty)$. Then
$$\lim_{n\rightarrow \infty}\frac{G_{n,A}(\Gamma)}{n}=\lim_{n\rightarrow \infty}\frac{G'_{n,A}(\Gamma)}{n}=d_{\Gamma, A}$$
for $c_{\Gamma, A}$ a constant that depends on $\Gamma$ and $A$.
\end{proposition}

In \cite{MR2770552},  Kahle defined and studied two models of random simplicial complexes over $G(n,r)$, the random Vietoris-Rips complex $R(n,r)$  and the random \v{C}ech complex $C(n,r)$. We define these random simplicial complexes  now. 

\begin{definition}
  The random \v{C}ech complex $C(n,r)$, also written
  $C(\mathcal{X}_n,r)$, is the simplicial complex with vertex set
  $\mathcal{X}_n$, and a  $k$-simplex $\sigma=(X_0,\dots,X_{k})\in C(n,r)$ iff
$$\bigcap_{X_i\in \sigma} B(X_i, \frac{r}{2})\neq \emptyset.$$
\end{definition}

\begin{definition}
The random Vietoris-Rips complex $R(n,r)$, also written $R(\mathcal{X}_n,r)$, is the simplicial complex with vertex set $\mathcal{X}_n$, and a $k$-simplex $\sigma=(X_0,\dots,X_{k})\in C(n,r)$ iff
$$B(X_i,\frac{r}{2})\cap B(X_j,\frac{r}{2})\neq \emptyset$$
for every pair $X_i,X_j\in \sigma.$
\end{definition}

The main idea of this article is to deduce similar results for the soft random simplicial complexes as the ones proved for $C(n,r)$ and $R(n,r)$ by Kahle in \cite{MR2770552}. We now present a  brief summary of the main arguments used by Kahle\cite{MR2770552}, which we will adapt to soft random graphs in this article. 

\begin{definition}[Definition 3.3 \cite{MR2770552}]\label{defO_k}
The $(k+1)$-dimensional cross-polytope is defined to be the convex hull of the $2k+2$ points $\{\pm e_i\}$, where $e_1,\dots,e_{k+1}$ are the standard basis vectors of $\mathbb{R}^{k+1}$. The boundary of this polytope is a k-dimensional simplicial complexes denoted by $O_k$.
\end{definition}
\begin{proposition}\label{cuentaso_k} Let $O_k$ as before, then:
\begin{equation}\label{xx1}
f_i(O_k)=2^{i+1}\binom{k+1}{i+1}.
\end{equation}
\end{proposition}

\begin{proof}
We will use induction on $k$ to prove (\ref{xx1}).
First notice that, $f_0(O_1)=4$ and $f_1(O_1)=4$  correspond with the result of equation (\ref{xx1}) for k=1.  Assume now that:
$$f_i(O_{k-1})=2^{i+1}\binom{k}{i+1}.$$
Let $f_i^j$ be the number of $i$-faces in $O_j$. By definition $O_k$ can be obtained from $O_{k-1}$ by adding 2 new vertices and including all possible new faces formed. Thus, since each $i-1$-face can span a new $i$-faces using each of the new vertices: 
$$f_i(O_k)=f_i^k=f_i^{k-1}+2f_{i-1}^{k-1}$$
therefore
\begin{align*}
f_i^k&=2^{i+1}\binom{k}{i+1}+2\times 2^{i}\binom{k}{i}\\
f_i^k&=2^{i+1}\left(\binom{k}{i+1}+\binom{k}{i}\right)\\
f_i(O_k)&=2^{i+1}\binom{k+1}{i+1}    
\end{align*}
\end{proof}

In \cite{MR2770552}, Kahle exhibited a smart argument to deduce information for the Betti numbers. He provided precise and adequate bounds for the Betti numbers  and deduced the limit behavior of these topological invariants by leveraging the corresponding limit behavior of the established bounds. Kahle observed that the $k$-th Betti number can be lower-bounded by the count of components of $G(n,r)$ isomorphic to $O_k$, since it represents the minimal complex required to increase $\beta_k$ by one. Furthermore, any contribution not coming from a component of $G(n,r)$ isomorphic to $O_k$ it must  comes from a component of at least $3k+2$ vertices and at least one $k$-face. Thus, if  $\hat{o}_k$ is the number of components isomorphic to $O_k$, and $f_k^{\geq 2k+3}$ is the number of $k$ faces in components of at least $2k+3$ vertices:
\begin{equation}\label{ineq1}
\hat{o}_k\leq \beta_k(R(n,r))\leq \hat{o}_k+ f_k^{\geq 2k+3},
\end{equation}
and therefore
$$E(\hat{o}_k)\leq E(\beta_k(R(n,r)))\leq E(\hat{o}_k)+ E(f_k^{\geq 2k+3}),$$
where  $\hat{o}_k$ is the number of components isomorphic to $O_k$, and $f_k^{\geq 2k+3}$ is the number of $k$ faces in components of at least $2k+3$ vertices. This inequality can be found in \cite{MR2770552} as equation (3.1). The inequalities 
\begin{equation}\label{ineq2}
\hat{S}_{k+1}\leq \beta_k(C(n,r))\leq\hat{S}_{k+1}+ f_k^{\geq k+3}
\end{equation}
and
$$E(\hat{S}_{k+1})\leq E(\beta_k(C(n,r)))\leq E(\hat{S}_{k+1})+ E(f_k^{\geq k+3}) $$
were proved similarly, and were used in \cite{MR2770552} to study the homology for $C(n,r)$. There, $\hat{S}_{k+1}$ denote the number of connected components isomorphic to the boundary of a $(k+1)$-simplicial complex and $f_k^{\geq k+3}$ is the number of $k$ faces in components of at least $k+3$ vertices. These inequalities, together with some of Penrose's results \cite{MR1986198}, were enough to obtain most of the results presented by Kahle in \cite{MR2770552}.

\begin{lemma}[Lemma 3.3 \cite{MR3079211} ]\label{sk} For $k\leq 1$, let $\hat{S}_{k+1}$ and  $S_{k+1}$ denote the number of connected components and subcomplexes isomorphic to the boundary of a $(k+1)$-simplicial complex and
$$\nu_{k+1}=\frac{1}{(k+1)!}\int_{\mathbb{R}^d} f(x)^k~dx \int_{(\mathbb{R}^d)^{k-1}}h_{1}(\{0,x_1,\dots,x_{k}\})~dx_1dx_2\dots dx_{k-1}$$
where $h_{1}(\{x_1,x_2,\dots,x_{k+1}\})$ is the characteristic function that is 1 when the set $\{x_1,x_2,\dots,x_{k+1}\}$ spans a simplicial complexes isomorphic to the boundary of a $(k+1)$-simplicial complex in $C(n,1)$. Then 
$$\lim_{n\rightarrow \infty}\frac{\hat{S}_{k+1}}{n^k r^{d(k-1)}}=\lim_{n\rightarrow \infty}\frac{S_{k+1}}{n^k r^{d(k-1)}}=\nu_{k+1}$$
\end{lemma}

The following theorem in the analogue of Lemma~\ref{sk} for the case of a general connected subcomplex X of $C(n,r)$.
\begin{theorem}\label{newcosito}
For $k\geq 1$, let $X$ a connected simplicial complex on $k$ vertices.  Let $h_X(x_1,...,x_{k})$ the characteristic function that is 1 if the set $\{x_1,...,x_{k+1}\}$ spans a subcomplex of $C(n,r)$ isomorphic to $X$ and 
$$\upsilon_X=\int_{\mathbb{R}^d}f(x)^kdx \int_{\mathbb{R}^d} h_X(0,y_1,...,y_{k}) dy_1...dy_{k-1}.$$
If $G_{n,r}(X)$ and $J_{n,r}(X)$ are, respectively, the number of subcomplex and components of $C(n,r)$ isomorphic to $X$, then  
$$\lim_{n\to \infty}\frac{E(G_{n,r}(X))}{n^{k}r^{d(k-1)}}=\lim_{n\to \infty}\frac{E(J_{n,r}(X))}{n^{k}r^{d(k-1)}}=\frac{\upsilon_X}{k!}.$$
\end{theorem}
\begin{proof}
$$E(G_{n,r}(X))=\binom{n}{k}\int_{\mathbb{R^d}}...\int_{\mathbb{R^d}}h_{X}(x_1,...,x_{k})\prod_{i=1}^{k}f(x_i)dx_1...dx_{k}$$
$$=\int_{\mathbb{R^d}}...\int_{\mathbb{R^d}}h_{X}(x_1,...,x_{k})f(x_1)^{k} dx_1...dx_{k}$$
$$ +\int_{\mathbb{R^d}}...\int_{\mathbb{R^d}}h_{X}(x_1,...,x_{k})\left(\prod_{i=1}^{k}f(x_i)-f(x_1)^{k}\right)dx_1...dx_{k}.$$
Now, as in proof of Proposition 1.3 from \cite{MR1986198},
$$ \int_{\mathbb{R^d}}...\int_{\mathbb{R^d}}h_{X}(x_1,...,x_{k})\left(\prod_{i=1}^{k}f(x_i)-f(x_1)^{k}\right)dx_1...dx_{k}\to 0.$$
For the first part, making the change of variables $x_i=x_1+ry_i$ for $2\leq i\leq k$:
$$\int_{\mathbb{R^d}}...\int_{\mathbb{R^d}}h_{X}(x_1,...,x_{k})f(x_1)^{k} dx_1...dx_{k}$$
$$=r^{d(k-1)}\int_{\mathbb{R^d}}...\int_{\mathbb{R^d}}h_{X}(x_1,x_1+ry_2,...,x_1+ry_{k})f(x_1)^k dy_2...dy_{k}dx_1$$
again, as in prove proof of Proposition 1.3 from \cite{MR1986198}, the above is asymptotic to $\nu_{X}r^{d(k-1)}$. Thus 
$$E(G_{n,r}(X))\sim \binom{n}{k}\nu_{X}r^{d(k-1)} \sim \frac{n^k}{k!}\nu_{X}r^{d(k-1)}$$.
On the other hand, let $\hat{h}_{X}(x_1,...,x_{k})$ the characteristic function that is 1 if the set $\{x_1,...,x_{k+1}\}$ spans a component of $C(n,r)$ isomorphic to $X$. Then
$$E(J_{n,r}(X))=\binom{n}{k}P(\hat{h}_{X}(x_1,...,x_{k})=1).$$
Now, the conditional probability $Q$ of the set $\{x_1,...,x_{k}\}$ spanning a component isomorphic to $X$ in $C(n,r)$ given that spans a subcomplex isomorphic to $X$ satisfies that
$$Q=P(\hat{h}_{X}(x_1,...,x_{k})=1\mid h_{X}(x_1,...,x_{k})=1 )\to 1$$
since 
$$Q\geq P(\text{no vertex $x_j \notin \{x_1,...,x_{k}\}$ is connected with $\{x_1,...,x_{k+1}\}$} )$$
on the other hand, the probability $Q_j$ that a vertex $x_j \notin \{x_1,...,x_{k}\}$ is connected with a vertex in $\{x_1,...,x_{k}\}$. Thus
$$Q_j=\int_{\bigcup\limits_{i=1}^kB(x_i,r)}f(x_j)dx_j\leq ||f||_{\infty}\int_{B(x_1,kr)}dx_j \leq ||f||_{\infty}(kr)^d\int_{B(0,1)}dx_j $$
since $X$ is connected. Therefore, if $\theta_d=Leb(B(0,1))$,  by independence 
$$Q\geq (1-Q_j)^{n-k}\geq (1-||f||_{\infty}(kr)^d \theta_d)^{n-k}\to 1.$$
Thus
\begin{align*}
E(J_{n,r}(X))&=\binom{n}{k}P(\hat{h}_{X}(x_1,...,x_{k})=1)\\
&=\binom{n}{k}P(h_{X}(x_1,...,x_{k})=1)Q\\
&=\binom{n}{k}P(h_{X}(x_1,...,x_{k})=1)\\
&= E(G_{n,r}(X))
\end{align*}

\end{proof}

\subsection{Main Results}

In this section, we provide a brief summary of the results found by Kahle for $R(n,r)$ and $C(n,r)$ in \cite{MR2770552}. Additionally, we present our main results, which focus on how the parameter $\rho$ must behave in order to maintain the results found by Kahle.
\begin{enumerate}
\item \textbf{The subcritical regime} The subcritical regime, also known as the sparse regime, is when $nr^d\to 0$. In this regime, the homology of the simplicial complexes $R(n,r)$ and $C(n,r)$ appears for the first time \cite{MR2770552}. With an additional hypothesis for the multiparameter $\rho$, this is the also the case for the homology of the soft random simplicial complexes $R(n,r,\rho)$ and $C(n,r,\rho)$. The following results shows the threshold radius for  the phase transition where homology appears for the soft random simplicial complexes.

\begin{stheorem}[\ref{theo6.4}]
Let $R(n,r,\rho)$ the  soft random  Vietoris-Rips complex with respect to a multiparameter vector $\rho=(p_1,p_2,\dots)$. Assume that $\lim\limits_{n\to \infty} \frac{nr^d}{ (1-p_{k+1})\prod\limits_{i=1}^{k}p_i^{\binom{k+1}{i}}}=0$.
\begin{enumerate}
\item If   $(1-p_{k+1})\prod\limits_{i=1}^{k}p_i^{\binom{k+2}{i+1}}n^{k+2}r^{d{k+1}}\to 0$, then $\beta_k(R(n,r,\rho))=0$ a.a.s
\item If  $(1-p_{k+1})\prod\limits_{i=1}^{k}p_i^{\binom{k+2}{i+1}}n^{k+2}r^{d{k+1}}\to \infty$, then $\beta_k(R(n,r,\rho))\neq 0$ a.a.s
\end{enumerate}
\end{stheorem}

\begin{stheorem}[\ref{theo6.10}]
Let $C(n,r,\rho)$ the  soft random \v{C}ech complex with respect to a multiparameter vector $\rho=(p_1,p_2,\dots)$. Assume that $\lim\limits_{n\to \infty} \frac{nr^d}{ \prod_{i=1}^{k}p_i^{\binom{k+1}{i}} }=0$.
\begin{enumerate}
\item If   $\prod\limits_{i=1}^{k}p_i^{\binom{k+2}{i+1}}n^{k+2}r^{d{k+1}}\to 0$, then $\beta_k(C(n,r,\rho))=0$ a.a.s
\item If  $\prod\limits_{i=1}^{k}p_i^{\binom{k+2}{i+1}}n^{k+2}r^{d{k+1}}\to \infty$, then $\beta_k(C(n,r,\rho))\neq 0$ a.a.s
\end{enumerate}
\end{stheorem}

\item \textbf{The critical regime} The critical regime, also known as the thermodynamic limit, is when $nr^d\to \lambda \in(0,\infty)$. In this regime, the number of components of the random simplicial complexes $R(n,r)$ and $C(n,r)$ is $\Theta(n)$ but significantly lower than in the subcritical regime \cite{MR2770552}. The following result ensure this is also the case for the soft random simplicial complexes.

\begin{stheorem}[\ref{theo4.1}] Consider the soft random simplicial complexes $R(n,r,\rho)$ and $C(n,r,\rho)$. Let $\epsilon>0$ be given. Assume that $nr^d \rightarrow \lambda\in (0,\infty)$. Then
\begin{enumerate}
\item  If $(1-p_{k+1})\prod\limits_{i=1}^{k}p_i^{\binom{k+2}{i+1}}\geq \frac{1}{n^{\epsilon}}$, 
then $E(\beta_k(R(n,r,\rho)))=O(n)$ and $E(\beta_k(R(n,r,\rho)))=\Omega(n^{1-\epsilon})$
\item If $\prod\limits_{i=1}^{k}p_i^{\binom{k+2}{i+1}}\geq \frac{1}{n^{\epsilon}}$, 
then $E(\beta_k(C(n,r,\rho)))=O(n)$ and $E(\beta_k(C(n,r,\rho)))=\Omega(n^{1-\epsilon})$
\end{enumerate}
\end{stheorem}

\item \textbf{The supercritical regime} The supercritical regime, sometimes called  the dense regime, is when $nr^d\to \infty$. In this regime, the number of components decay slowly \cite{MR2770552}. The following result ensure this is the case too for the soft  random Vietoris-Rips  complex.

\begin{stheorem}[\ref{theo5.1}]
Let $R(n,r,\rho)$ be the soft random  Vietoris-Rips complex for which the  $n$ points $ \{X_1,\dots,X_n\}$ are taken i.i.d uniformly from a  bounded convex body $K \subset \mathbb{R}^d$. Suppose $r=\omega(n^{-\frac{1}{d}})$ and write $W=nr^d$. Then
$$E( \beta_k(R(n,r,\rho)))=O\left(nW^k\prod_{i=1}^{k}p^{\binom{k+1}{i+1}}e^{-cW\prod\limits_{i=1}^{k+1}p^{\binom{k+1}{i}}}\right).$$
\end{stheorem}

\item \textbf{The connected regime} The connected regime is when $nr^d$ grows faster than $\log(n)$. It is contained in the supercritical regime, and it is the regime where the random  simplicial complexes $R(n,r)$ and $C(n,r)$  become highly connected. Therefore, $R(n,r)$ and $C(n,r)$ exhibit a second threshold where their $k$-homology changes from non-zero to zero  \cite{MR2770552}. The following results confirm that this is also true for the soft random simplicial complexes. 

\begin{stheorem}[\ref{theo666}]
Let $R(n,r,\rho)$ be the soft Vietoris-rips for which the  $n$ points were taken i.i.d uniformly from a  bounded convex body $K$ in $\mathbb{R}^d$. Fix $k\geq 0$, let $r\geq c_k(\frac{\log(n)}{n})^{\frac{1}{d}}$ , and suppose that for each $i\leq k$, $p_i(n)$ is a sequence bounded away from zero. Then $R(n,r,\rho)$ is $k$-connected a.a.s.
\end{stheorem}

\begin{stheorem}[\ref{theo5.6}]
Let $\epsilon>0$ be given and let $C(n,r,\rho)$ be the soft random \v{C}ech complex for which  $n$ points were taken i.i.d uniformly from a  bounded convex body $K$ in $\mathbb{R}^d$.  If $p_i\geq 1-\frac{1}{n^{2+i+\epsilon}}$ for all $i$, then there exists a constant $c_K$ depending on $K$, such that if $r\geq c_K(\frac{\log(n)}{n})^{\frac{1}{d}}$, then $C(n,r,\rho)$ is contractible a.a.s.
\end{stheorem}
The following results summarize the changes in homology observed in the soft random simplicial complexes discussed in this document.
\begin{theorem}[\ref{need1}]
Let $R(n,r,\rho)$ be the soft Vietoris-rips complex for which  the  $n$ points were taken i.i.d uniformly from a  bounded convex body $K \subset \mathbb{R}^d$. For any fixed $k\geq 1$:
\begin{enumerate}
\item If either $r^{d{k+1}}=o((1-p_{k+1})^{-1}\prod\limits_{i=1}^{k}p_i^{-\binom{k+2}{i+1}}n^{-(k+2)})$ or else $r= \omega((\frac{\log(n)}{n})^{\frac{1}{d}})$ and  for all $i\leq k$, $p_i(n)$ is a sequence bounded away from zero, then $H_k=0$ a.a.s.
\item If $r^{d{k+1}}=\omega((1-p_{k+1})^{-1}\prod\limits_{i=1}^{k}p_i^{-\binom{k+2}{i+1}}n^{-(k+2)})$, $r= o((\frac{\log(n)}{n})^{\frac{1}{d}})$, and for all $i\leq k$, $p_i(n)\to 1$, then $H_k\neq 0$ a.a.s.
\end{enumerate}
\end{theorem}

\begin{stheorem}[\ref{need2}]
Let $C(n,r,\rho)$ be the soft \v{C}ech complex for which the  $n$ points were taken i.i.d uniformly from a bounded convex body $K \in \mathbb{R}^d$. For any fixed $k\geq 1$
\begin{enumerate}
\item If either $r^{d{k+1}}=o(\prod\limits_{i=1}^{k}p_i^{-\binom{k+2}{i+1}}n^{-(k+2)})$ or else $r= \omega((\frac{\log(n)}{n})^{\frac{1}{d}})$ and for all $i$, $p_i\geq 1-\frac{1}{n^{2+i+\epsilon}}$, then $H_k=0$ a.a.s.
\item If $r^{d{k+1}}=\omega(\prod\limits_{i=1}^{k}p_i^{-\binom{k+2}{i+1}}n^{-(k+2)})$, $r= o((\frac{\log(n)}{n})^{\frac{1}{d}})$, and for all $i\leq k$, $p_i(n)\to 1$, then $H_k\neq 0$  a.a.s.
\end{enumerate}
\end{stheorem}
\end{enumerate}

In this article, we introduce random variables associated to different versions of random simplicial complexes.  Thus, we need some notation that allows us to differentiate whether we are talking about a random variable related to the Poisson version of a random simplicial complex or not; or if we are talking about a random variable related to a random simplicial complex that depends on a probability vector $\rho=(p_1,p_2,\dots)$ i.e. having a probability $p_i$ that adds randomness to the corresponding dimension. Therefore, we make the following notational conventions: All random variables related to a  random simplicial complex that depends on a probability vector $\rho=(p_1,p_2,\dots)$ explicitly depend on this parameter. Also, all random variables related to  random graphs or random simplicial complexes built over the  Poisson process $\mathcal{P}_n$ will be indicated by the super index $\mathcal{P}$. Thus:

\begin{enumerate}
\item Random variables without an explicit dependence of the parameter $\rho$ and without the super index $\mathcal{P}$ are associated to the graph $G(n,r)$ or to the simplicial complexes $R(n,r)$ and $C(n,r)$ i.e. to random graphs or random complexes that do not depend on a probability vector $\rho$ and whose vertex set are $n$ independent distributed points  $X_1,\dots,X_n$.
\item Random variables without an explicit dependence of the parameter $\rho$ and with the super index $\mathcal{P}$ are associated to the graph $G(\mathcal{P}_n,r)$ or to the simplicial complexes $R(\mathcal{P}_n,r)$ and $C(\mathcal{P}_n,r)$ i.e to random graphs or random complexes that do not depend on a probability vector $\rho$ and whose vertices depend on a Poisson process $\mathcal{P}_n$.
\item Random variables with an explicit dependence of the parameter $\rho$ and without the super index $\mathcal{P}$ are associated to $R(n,r,\rho)$ or to $C(n,r,\rho)$ i.e to  random complexes that do depend on a probability vector $\rho$ and whose vertex set are $n$ independent distributed points  $X_1,\dots,X_n$.
\item Random variables with an explicit dependency of the parameter $\rho$ and with the super index $\mathcal{P}$ are associated to $R(\mathcal{P}_n,r,\rho)$ or to $C(\mathcal{P}_n,r,\rho)$. i.e. to  random complexes that do depend on a probability vector $\rho$ and whose vertex set depends on a  Poisson process $\mathcal{P}_n$.
\end{enumerate}

Finally, we use Bachmann-Landau Big-O, little-o notation. In particular, we say that for two non-negative functions $f(n)$ and $g(n)$:
\begin{itemize}
\item $g(n)=\Omega(f(n))$, iff there exists $N$ and $k$ such that for all $n>N$, $ k f(n)\leq g(n)$.
\item $g(n)=\Theta(f(n))$, iff there exists $N$, $k_1$ and $k_2$ such that for all $n>N$, $ k_1 f(n)\leq g(n)\leq k_2 f(n)$.
\item $g(n)=o(f(n))$, iff for all $\epsilon>0$, there exists $N$ such that for all $n>N$, $g(n)\leq \epsilon f(n)$.
\item $g(n)=\omega(f(n))$, iff for all $\epsilon>0$, there exists $N$ such that for all $n>N$, $ \epsilon f(n)\leq g(n)$.
\item We also say that $f(n)\sim g(n)$ iff
$$\lim_{n\to\infty}\frac{f(n)}{(g(n))}=1.$$
\end{itemize}

\section{Soft random simplicial complexes}\label{section2}
In this section, we define the soft random simplicial complexes $R(n,r,\rho)$ and $C(n,r,\rho)$, and we establish the results we will use in order to study them.
\begin{definition}
Let $\rho$ be the infinite multiparameter vector $\rho=(p_1,p_2,\dots)$ with $0\leq p_i\leq 1$ for each $i$. 
\begin{enumerate}
\item The soft random Vietoris-Rips complexes $R(n,r,\rho)$  is the simplicial complex with vertex set $\mathcal{X}_n$ in which a $k$-simplex $\sigma=(X_0,\dots,X_{k})$ in included on it with probability $p_k$ iff:
$$B(X_i,\frac{r}{2})\cap B(X_j,\frac{r}{2})\neq \emptyset$$
for every pair $X_i,X_j\in \sigma$.
\item The soft random \v{C}ech complex $C(n,r,\rho)$ is the simplicial complex with vertex set $\mathcal{X}_n$ in which a $k$-simplex $\sigma=(X_0,\dots,X_{k})$ in included in it with probability $p_k$ iff:
$$\bigcap_{X_i\in \sigma} B(X_i, \frac{r}{2})\neq \emptyset.$$
\end{enumerate}
\end{definition}

The following  results  relate the expected value of the number of subcomplexes of $R(n,r,\rho)$ isomorphic to a given connected simplicial complex $Y$ and the expected value of the number of subcomplexes isomorphic to $Y$ in  the original random geometric complex $R(n,r)$.

\begin{theorem}\label{theo6.1}
Let $\rho$ be a multiparameter vector $\rho=(p_1,p_2,\dots)$ and  let $R(n,r,\rho)$ be the soft random Vietoris-Rips complex. Also, assume that $Y$ is a connected simplicial complex on $k$ vertices that is realizable in both $R(n,r,\rho)$  and $R(n,r)$ and $G_{n,\rho}(Y)$ is the number of subcomplex of  $R(n,r,\rho)$  isomorphic to $Y$. Then
$$ c_1\prod_{i=1}^{k-1}p_i^{f_i(Y)} E(G_{n}(Y))\leq E(G_{n,\rho}(Y))\leq c_2 \prod_{i=1}^{k-1}p_i^{f_i(Y)}E(G_{n}(Y))$$
where $G_{n}(Y)$ is the number of subcomplex of  $R(n,r)$  isomorphic to $Y$ and $c_1,c_2$ are constants depending on $k$.
\end{theorem}

\begin{proof}
Let $Y$ be the subcomplex of $R(n,r,\rho)$ on the vertices $\{v_1,\dots,v_k\}$. Let $\mathcal{A}$ be the collection of all the simplicial complexes on vertices $\{v_1,\dots,v_k\}$ that contain $Y$. For $\hat{Y}$ a simplicial complex on $k$ vertices, define the characteristic functions $H_{\hat{Y}}$ and $Q_{\hat{Y}}$ as follows
\begin{enumerate}
\item $H_{\hat{Y}}(v_1,\dots,v_k)=1$ iff $\hat{Y}$ is the subcomplex spanned by $\{v_1,\dots,v_k\}$ on  $R(n,r)$.
\item $Q_{\hat{Y}}(v_1,\dots,v_k)=1$ iff $\hat{Y}$ is the subcomplex spanned by $\{v_1,\dots,v_k\}$ on  $R(n,r,\rho)$.
\end{enumerate}
Given a simplicial complex $\hat{Y}\in \mathcal{A}$, consider the probability of  $Y$ appearing in $R(n,r,\rho)$ given that  $\hat{Y}$ is in $R(n,r)$:
$$P(Q_{Y}(v_1,\dots,v_k)=1| H_{\hat{Y}}(v_1,\dots,v_k)=1)=\prod_{i=1}^{k-1}p_i^{f_i(Y)}(1-p_i)^{e_i(\hat{Y},Y)}$$
where $e_i(\hat{Y},Y)$ denotes the number of $i$-faces of $\hat{Y}$ whose $(i-1)$-skeletons are contained in $Y$. Therefore summing over all possible $\hat{Y}$ which contain $Y$, we have:
$$P(Q_{Y}(v_1,\dots,v_k)=1)=\sum_{Y\subset\hat{Y}}\prod_{i=1}^{k-1}p_i^{f_i(Y)}(1-p_i)^{e_i(\hat{Y},Y)}P(H_{\hat{Y}}(v_1,\dots,v_k)=1)$$
then, since $(1-p_i)\leq 1$:
$$P(Q_{Y}(v_1,\dots,v_k)=1)\leq\sum_{Y\subset\hat{Y}}\prod_{i=1}^{k-1}p_i^{f_i(Y)}P(H_{\hat{Y}}(v_1,\dots,v_k)=1)$$
now, in $R(n,r)$ every subsimplicial complex $\hat{Y}$  is totally determined by its $1$-skeleton, thus if $\hat{\Gamma}$ is the graph obtained by the edges and the vertices of $\hat{Y}$, $P(H_{\hat{Y}}(v_1,\dots,v_k)=1)=P(\hat{\Gamma})$ . Therefore,  if $\Gamma$ is the graph obtained from $Y$
$$P(Q_{Y}(v_1,\dots,v_k)=1)\leq \sum_{\Gamma\subset\hat{\Gamma}}\prod_{i=1}^{k-1}p_i^{f_i(Y)}P(\hat{\Gamma}).$$
By Proposition~\ref{Penrose-1}, $P(\hat{\Gamma})=c_{\hat{\Gamma}}r^{d(k-1)}$ where $c_{\hat{\Gamma}}$ is a constant depending on $\hat{\Gamma}$. Since the number of possible graphs $\hat{\Gamma}$ that contain $\Gamma$ is finite:
$$P(Q_{Y}(v_1,\dots,v_k)=1)\leq c \prod_{i=1}^{k-1}p_i^{f_i(Y)}r^{d(k-1)}$$
thus
$$\binom{n}{k}P(Q_{Y}(v_1,\dots,v_k)=1)\leq c \binom{n}{k} \prod_{i=1}^{k-1}p_i^{f_i(Y)}r^{d(k-1)}$$
therefore, again by  Proposition~\ref{Penrose-1}:
$$E(G_{n,\rho}(Y))\leq c_2 \binom{n}{k}\prod_{i=1}^{k-1}p_i^{f_i(Y)}r^{d(k-1)}\sim c_2\prod_{i=1}^{k-1}p_i^{f_i(Y)} E(G_n(Y))$$
Finally, for the lower bound since
$$P(Q_{Y}(v_1,\dots,v_k)=1)=\sum_{Y\subset\hat{Y}}\prod_{i=1}^{k-1}p_i^{f_i(Y)}(1-p_i)^{e_i(\hat{Y},Y)}P(H_{\hat{Y}}(v_1,\dots,v_k)=1)$$
then
\begin{align*}
P(Q_{Y}(v_1,\dots,v_k)=1) &= \prod_{i=1}^{k-1}p_i^{f_i(Y)}P(H_{Y}(v_1,\dots,v_k)=1) \\
&\quad + \sum_{Y\subset\hat{Y}\neq Y} \prod_{i=1}^{k-1}p_i^{f_i(Y)}(1-p_i)^{e_i(\hat{Y},Y)} P(H_{\hat{Y}}(v_1,\dots,v_k)=1)
\end{align*}
therefore
\begin{align*}
P(Q_{Y}(v_1,\dots,v_k)=1) &\geq  \prod_{i=1}^{k-1}p_i^{f_i(Y)}P(H_{Y}(v_1,\dots,v_k)=1) 
\end{align*}
thus
$$P(Q_{Y}(v_1,\dots,v_k)=1)\geq \prod_{i=1}^{k-1}p_i^{f_i(Y)}P(\Gamma)$$
therefore, asymptotically by Proposition~\ref{Penrose-1}:
$$E(E(G_{n,\rho}(Y)))\geq  c_1\prod_{i=1}^{k-1}p_i^{f_i(Y)} E(G_n(Y)).$$
\end{proof}

\begin{theorem}\label{theo6.2}
Let $\rho$ be a multiparameter vector $\rho=(p_1,p_2,\dots)$ and $R(n,r,\rho)$ the soft random simplicial Vietoris-Rips complex. Also, assume $Y$ is a connected simplicial complex on $k$ vertices that is realizable in both $R(n,r,\rho)$  and $R(n,r)$. Then, if $J_{n,\rho}(Y)$ is the number of subcomplex of  $R(n,r,\rho)$  isomorphic to $Y$: 
$$ \prod_{i=1}^{k-1}p_i^{f_i(Y)} E(J_n(Y))\leq E(J_{n,\rho}(Y))\leq c \prod_{i=1}^{k-1}p_i^{f_i(Y)}E(J_n(Y))$$
where $J_{n,\rho}(Y)$ is the number of subcomplex of  $R(n,r)$  isomorphic to $Y$ and $c_1,c_2$ are constants depending on $k$.
\end{theorem}

\begin{proof}
Let $Y$ be an isolated and connected subcomplex of $R(n,r,\rho)$ on the vertices $\{v_1,\dots,v_k\}$. Let $\mathcal{A}$ be the collection of all the simplicial complexes on vertices $\{v_1,\dots,v_k\}$ that contain $Y$. For $\hat{Y}$ a simplicial complex on $k$ vertices, define the characteristic functions $H_{\hat{Y}}$ and $Q_{\hat{Y}}$ as follows
\begin{enumerate}
\item $H_{\hat{Y}}(v_1,\dots,v_k)=1$ iff $\hat{Y}$ is a component spanned by $\{v_1,\dots,v_k\}$ on  $R(n,r)$.
\item $Q_{\hat{Y}}(v_1,\dots,v_k)=1$ iff $\hat{Y}$ is a component spanned by $\{v_1,\dots,v_k\}$ on  $R(n,r,\rho)$.
\end{enumerate}
Given a simplicial complex $\hat{Y}\in \mathcal{A}$, consider the probability of  $Y$ appearing in $R(n,r,\rho)$ given that  $\hat{Y}$ is in $R(n,r)$:
$$P\left(Q_{Y}(v_1,\dots,v_k)=1| H_{\hat{Y}}(v_1,\dots,v_k)=1\right)=\prod_{i=1}^{k-1}p_i^{f_i(Y)}(1-p_i)^{e_i(\hat{Y},Y)}$$
where $e_i(\hat{Y},Y)$ denotes the number of $i$-faces of $\hat{Y}$ whose $(i-1)$-skeletons are contained in $Y$. Therefore summing over all possible $\hat{Y}$ which contain $Y$, we have:
$$P\left(Q_{Y}(v_1,\dots,v_k)=1\right)=\sum_{Y\subset\hat{Y}}\prod_{i=1}^{k-1}p_i^{f_i(Y)}(1-p_i)^{e_i(\hat{Y},Y)}P(H_{\hat{Y}}(v_1,\dots,v_k)=1)$$
then, since $(1-p_i)\leq 1$:
$$P\left(Q_{Y}(v_1,\dots,v_k)=1\right)\leq\sum_{Y\subset\hat{Y}}\prod_{i=1}^{k-1}p_i^{f_i(Y)}P(H_{\hat{Y}}(v_1,\dots,v_k)=1)$$
now, in $R(n,r)$ every subsimplicial complex $\hat{Y}$  is totally determined by its $1$-skeleton, thus if $\hat{\Gamma}$ is the graph obtained by the edges and the vertices of $\hat{Y}$, $P(H_{\hat{Y}}(v_1,\dots,v_k)=1)=P(\hat{\Gamma})$ . Therefore,  if $\Gamma$ is the graph obtained from $Y$
$$P\left(Q_{Y}(v_1,\dots,v_k)=1\right)\leq\sum_{\Gamma\subset\hat{\Gamma}}\prod_{i=1}^{k-1}p_i^{f_i(Y)}P(\hat{\Gamma}).$$
By Proposition~\ref{Penrose-2}, $P(\hat{\Gamma})=c_{\hat{\Gamma}}r^{d(k-1)}$ where $c_{\hat{\Gamma}}$ is a constant depending on $\hat{\Gamma}$. Since the number of possible graphs $\hat{\Gamma}$ that contain $\Gamma$ is finite:
$$P\left(Q_{Y}(v_1,\dots,v_k)=1\right)\leq c \prod_{i=1}^{k-1}p_i^{f_i(Y)}r^{d(k-1)}$$
thus
$$\binom{n}{k}P\left(Q_{Y}(v_1,\dots,v_k)=1\right)\leq c \binom{n}{k} \prod_{i=1}^{k-1}p_i^{f_i(Y)}r^{d(k-1)}$$
therefore, again by  Proposition~\ref{Penrose-2}:
$$E(J_{n,\rho}(Y))\leq c_2 \binom{n}{k}\prod_{i=1}^{k-1}p_i^{f_i(Y)}r^{d(k-1)}\sim c_2\prod_{i=1}^{k-1}p_i^{f_i(Y)} E(J_n(Y))$$
Finally, for the lower bound since
$$P\left(Q_{Y}(v_1,\dots,v_k)=1\right)=\sum_{Y\subset\hat{Y}}\prod_{i=1}^{k-1}p_i^{f_i(Y)}(1-p_i)^{e_i(\hat{Y},Y)}P(H_{\hat{Y}}(v_1,\dots,v_k)=1)$$
then
\begin{align*}
P(Q_{Y}(v_1,\dots,v_k)=1) &= \prod_{i=1}^{k-1}p_i^{f_i(Y)}P(H_{Y}(v_1,\dots,v_k)=1) \\
&\quad + \sum_{Y\subset\hat{Y}\neq Y} \prod_{i=1}^{k-1}p_i^{f_i(Y)}(1-p_i)^{e_i(\hat{Y},Y)} P(H_{\hat{Y}}(v_1,\dots,v_k)=1)
\end{align*}
therefore
\begin{align*}
P(Q_{Y}(v_1,\dots,v_k)=1) &\geq \prod_{i=1}^{k-1}p_i^{f_i(Y)}P(H_{Y}(v_1,\dots,v_k)=1)
\end{align*}
thus
$$P\left(Q_{Y}(v_1,\dots,v_k)=1\right)\geq \prod_{i=1}^{k-1}p_i^{f_i(Y)}P(\Gamma)$$
therefore, asymptotically by Proposition~\ref{Penrose-2}:
$$E(J_{n,\rho}(Y))\geq  c_1\prod_{i=1}^{k-1}p_i^{f_i(Y)} E(J_n(Y)).$$
\end{proof}

\begin{theorem}\label{newadd}
Let $\rho$ be the multiparameter vector $\rho=(p_1,p_2,\dots)$ and let $R(n,r,\rho)$ be the corresponding soft random simplicial Vietoris-Rips complex. If $\hat{T}_{k+1}(\rho)$ and $\hat{T}_{k+1}$ denote the number of connected components in $R(n,r,\rho)$ and $R(n,r)$, respectively, isomorphic to the boundary of a $(k+1)$-simplex and isomorphic to a $(k+1)$-simplex. Then 
$$  E(\hat{T}_{k+1}(\rho))= (1-p_{k+1}) \prod_{i=1}^{k+1}p_i^{\binom{k+2}{i+1}}E(\hat{T}_{k+1}).$$
\end{theorem}

\begin{proof}
Let $\sigma$ be a simplicial complex on the vertices $\{v_1,\dots,v_{k+2}\}$ isomorphic to the boundary of a $k+1$ simplex. Let $\mathcal{A}$ be the collection of all the simplicial complexes on vertices $\{v_1,\dots,v_{k+2}\}$ that contain $\sigma$. For $Y$ a simplicial complex on $k+2$ vertices, define the characteristic functions $H_{\hat{Y}}$ and $Q_{\hat{Y}}$ as follows
\begin{enumerate}
\item $H_{Y}(v_1,\dots,v_{k+2})=1$ iff $Y$ is the subcomplex spanned by $\{v_1,\dots,v_{k+2}\}$ on  $R(n,r)$.
\item $Q_{Y}(v_1,\dots,v_{k+2})=1$ iff $Y$ is the subcomplex spanned by $\{v_1,\dots,v_{k+2}\}$ on  $R(n,r,\rho)$.
\end{enumerate}
Given a simplicial complex $Y\in \mathcal{A}$, consider the probability of  $\sigma$ appearing in $R(n,r,\rho)$ given that  $Y$ is in $R(n,r)$:

$$P\left(Q_{\sigma}(v_1,\dots,v_{k+2}\right)=1| H_{Y}(v_1,\dots,v_{k+2})=1)=\prod_{i=1}^{k}p_i^{f_i(\sigma)}(1-p_i)^{e_i(Y,\sigma)}$$
where $e_i(Y,\sigma)$ denotes the number of $i$-faces of $Y$ whose $(i-1)$-skeletons are contained in $\sigma$.Therefore
$$P\left(Q_{\sigma}(v_1,\dots,v_{k+2})=1\right)=\sum_{\sigma\subset Y}\prod_{i=1}^{k}p_i^{f_i(\sigma)}(1-p_i)^{e_i(Y,\sigma)}P(H_{Y}(v_1,\dots,v_{k+2})=1)$$
now, notice that in $R(n,r)$ there is .only one simplicial complex that properly contains $\sigma$, the simplex $\Delta$ on $k+2$ vertices,
so
$$P\left(Q_{\sigma}(v_1,\dots,v_{k+2})=1\right)=\prod_{i=1}^{k}p_i^{\binom{k+2}{i+1}}(1-p_{k+1})P(H_{\Delta}(v_1,\dots,v_{k+2})=1)$$
thus
\begin{align*}
E(\hat{T}_{k+1}(\rho))&=\binom{n}{k+2}(1-p_{k+1})\prod_{i=1}^{k}p_i^{\binom{k+2}{i+1}}P(H_{\Delta}(v_1,\dots,v_{k+2})=1).
\end{align*}
By Proposition~\ref{Penrose-1}, $P(H_{\Delta}(v_1,\dots,v_{k+2})=1)=c_{\Delta}r^{d(k+1)}$, therefore
$$E(\hat{T}_{k+1}(\rho))=\binom{n}{k+2}(1-p_{k+1})\prod_{i=1}^{k}p_i^{\binom{k+2}{i+1}} c_{\Delta} r^{d(k+1)} $$
then
$$E(\hat{T}_{k+1}(\rho))=(1-p_{k+1})\prod_{i=1}^{k}p_i^{\binom{k+2}{i+1}} c_{\Delta} r^{d(k+1)} $$
as desired, and the proof is complete.
\end{proof}

Now, we want to deduce similar results for the soft random \v{C}ech complex $C(n,r,\rho)$.  Notice that, we need to be very careful since in this case, the simplicial complex is not determined by its 1-skeleton. However, we can still use some arguments to deduce the theorems.

\begin{theorem}\label{theo6.6}
Let $\rho$ be the multiparameter vector $\rho=(p_1,p_2,\dots)$ and let $C(n,r,\rho)$ be the corresponding soft random simplicial \v{C}ech complex. If $\hat{S}_{k+1}(\rho)$ and $\hat{S}_{k+1}$ denote the number of connected components in $C(n,r,\rho)$ and $C(n,r)$, respectively, isomorphic to the boundary of a $(k+1)$-simplex. Then for some constants $c_1$ and $c_2$, depending on $k+2$, we have
$$ c_1\prod_{i=1}^{k}p_i^{\binom{k+2}{i+1}} E(\hat{S}_{k+1})\leq E(\hat{S}_{k+1}(\rho))\leq c_2 \prod_{i=1}^{k+1}p_i^{\binom{k+2}{i+1}}E(\hat{S}_{k+1}).$$
\end{theorem}

\begin{proof}
Let $\sigma$ be a simplicial complex on the vertices $\{v_1,\dots,v_{k+2}\}$ isomorphic to the boundary of a $k+1$ simplex. Let $\mathcal{A}$ be the collection of all the simplicial complexes on vertices $\{v_1,\dots,v_{k+2}\}$ that contain $\sigma$. For $Y$ a simplicial complex on $k+2$ vertices, define the characteristic functions $H_{\hat{Y}}$ and $Q_{\hat{Y}}$ as follows
\begin{enumerate}
\item $H_{Y}(v_1,\dots,v_{k+2})=1$ iff $Y$ is the subcomplex spanned by $\{v_1,\dots,v_{k+2}\}$ on  $C(n,r)$.
\item $Q_{Y}(v_1,\dots,v_{k+2})=1$ iff $Y$ is the subcomplex spanned by $\{v_1,\dots,v_{k+2}\}$ on  $C(n,r,\rho)$.
\end{enumerate}
Given a simplicial complex $Y\in \mathcal{A}$, consider the probability of  $\sigma$ appearing in $C(n,r,\rho)$ given that  $Y$ is in $C(n,r)$:

$$P\left(Q_{\sigma}(v_1,\dots,v_{k+2}\right)=1| H_{Y}(v_1,\dots,v_{k+2})=1)=\prod_{i=1}^{k}p_i^{f_i(\sigma)}(1-p_i)^{e_i(Y,\sigma)}$$
where $e_i(Y,\sigma)$ denotes the number of $i$-faces of $Y$ whose $(i-1)$-skeletons are contained in $\sigma$.Therefore
$$P\left(Q_{\sigma}(v_1,\dots,v_{k+2})=1\right)=\sum_{\sigma\subset Y}\prod_{i=1}^{k}p_i^{f_i(\sigma)}(1-p_i)^{e_i(Y,\sigma)}P(H_{Y}(v_1,\dots,v_{k+2})=1)$$
now, notice that in $C(n,r)$ there are only two simplicial complexes that contains $\sigma$, the simplex $\Delta$ on $k+2$ vertices  and $\sigma$ itself,
so
\begin{align}
P\left(Q_{\sigma}(v_1,\dots,v_{k+2})=1\right) 
&= \prod_{i=1}^{k}p_i^{\binom{k+2}{i+1}}P\left(H_{\sigma}(v_1,\dots,v_{k+2})=1\right) \nonumber \\
&\quad + \prod_{i=1}^{k}p_i^{\binom{k+2}{i+1}}(1-p_{k+1})P\left(H_{\Delta}(v_1,\dots,v_{k+2})=1\right)
\end{align}
thus
\begin{align*}
E(\hat{S}_{k+1}(\rho))&= \binom{n}{k+2}\prod_{i=1}^{k}p_i^{\binom{k+2}{i+1}}P(H_{\sigma}(v_1,\dots,v_{k+2})=1)\\
&+\binom{n}{k+2}\prod_{i=1}^{k}(1-p_{k+1})p_i^{\binom{k+2}{i+1}}P(H_{\Delta}(v_1,\dots,v_{k+2})=1).
\end{align*}
By Lemma~\ref{newcosito}, $P(H_{\Delta}(v_1,\dots,v_{k+2})=1)=c_{\Delta}r^{d(k+1)}$ and $P(H_{\sigma}(v_1,\dots,v_{k+2})=1)=c_{\sigma}r^{d(k+1)}$, therefore
$$E(\hat{S}_{k+1}(\rho))\leq \binom{n}{k+2}\prod_{i=1}^{k}p_i^{\binom{k+2}{i+1}}c_{\sigma} r^{d(k+1)}+\binom{n}{k+2}\prod_{i=1}^{k}p_i^{\binom{k+2}{i+1}} c_{\Delta} r^{d(k+1)} $$
thus, there exists a constant $c$ such that 
$$E(\hat{S}_{k+1}(\rho))\leq c r^{d(k+1)}\binom{n}{k+2}\prod_{i=1}^{k}p_i^{\binom{k+2}{i+1}}$$
also, by Lemma~\ref{sk} there exists a constant $c_2$ such that:
$$E(\hat{S}_{k+1}(\rho))\leq c r^{d(k+1)}\binom{n}{k+2}\prod_{i=1}^{k}p_i^{\binom{k+2}{i+1}}\sim c_2\prod_{i=1}^{k}p_i^{\binom{k+2}{i+1}}E(\hat{S}_{k+1}).$$
On the other hand if we ignore the second term in the sum we obtain 
$$E(\hat{S}_{k+1}(\rho))\geq c_2\prod_{i=1}^{k}p_i^{\binom{k+2}{i+1}}E(\hat{S}_{k+1})$$
as desired, and the proof is complete.
\end{proof}
\begin{theorem}\label{theo6.8}
Let $\sigma$ be a simplicial complex that is the connected union of a $k$-simplex and 2 edges and $s_j$ the number of subcomplex of $C(n,r,\rho)$ isomorphic to $\sigma$. Then for a constant $c$ depending on $k+3$
$$E(s_j)\leq cp_1^2\prod_{i=1}^{k+1}p_i^{\binom{k+1}{i+1}}E(s_j)$$
\end{theorem}
\begin{proof}
Let $\hat{\sigma}$ a $k$-simplex on the vertices $\{v_1,\dots,v_{k+1}\}$ and let $\sigma$ be the connected union of $\hat{\sigma}$ together with the two new edges and the vertices $v_{k+3}$ and $v_{k+4}$. Let $\mathcal{A}$ be the collection of all the simplicial complexes on vertices $\{v_1,\dots,v_{k+3}\}$ that contain $\sigma$. For $Y$ a simplicial complex on $k+3$ vertices, define the characteristic functions $H_{Y}$ and $Q_{Y}$ as follows
\begin{enumerate}
\item $H_{Y}(v_1,\dots,v_{k+3})=1$ iff $Y$ is the subcomplex spanned by $\{v_1,\dots,v_{k+3}\}$ on  $C(n,r)$.
\item $Q_{Y}(v_1,\dots,v_{k+3})=1$ iff $Y$ is the subcomplex spanned by $\{v_1,\dots,v_{k+3}\}$ on  $C(n,r,\rho)$.
\end{enumerate}
Given a simplicial complex $Y\in \mathcal{A}$, consider the probability of  $\sigma$ appearing in $C(n,r,\rho)$ given that  $Y$ is in $C(n,r)$:

$$P\left(Q_{\sigma}(v_1,\dots,v_{k+3})=1|H_{Y}(v_1,\dots,v_{k+3})=1 \right)=\prod_{i=1}^{k}p_i^{f_i(\sigma)}(1-p_i)^{e_i(Y,\sigma)}$$
where $e_i(Y,\sigma)$ denotes the number of $i$-faces of $Y$ whose $(i-1)$-skeletons are contained in $\sigma$.Therefore
$$P\left(Q_{\sigma}(v_1,\dots,v_{k+3})=1\right)=\sum_{\sigma\subset Y}\prod_{i=1}^{k}p_i^{f_i(\sigma)}(1-p_i)^{e_i(Y,\sigma)}P(H_{Y}(v_1,\dots,v_{k+3})=1)$$
now, for Lemma~\ref{newcosito}, all simplicial complexes $Y$ have probability $c_Yr^{d(k+2)}$. Since there are at most a finite amount of such Y's:
$$P\left(Q_{\sigma}(v_1,\dots,v_{k+3})=1\right)\leq c_2\prod_{i=1}^{k}p_i^{f_i(\sigma)}r^{d(k+2)} $$
thus
$$P\left(Q_{\sigma}(v_1,\dots,v_{k+3})=1\right)\leq c_2p_1^2\prod_{i=1}^{k}p_i^{\binom{k+1}{i+1}}r^{d(k+2)} $$
therefore
\begin{equation}\label{ineq10}
E(s_j(\rho))\leq  c_2p_1^2 \binom{n}{k+3}\prod_{i=1}^{k}p_i^{\binom{k+1}{i+1}}r^{d(k+2)}.
\end{equation}
On the other hand,  $E(s_j)\sim n^{k+3}r^{d(k+2)}$. Thus, asymptotically
$$E(s_j(\rho))\leq  c_2p_1^2 \prod_{i=1}^{k}p_i^{\binom{k+1}{i+1}}E(s_j).$$
\end{proof}

\section{Subcritical regime}\label{section3}

Motivated by the way in which the geometric graph $G(n,r)$ and the simplicial complexes $C(n,r)$ and $R(n,r)$ presented different behaviors in the so-called subcritical, critical, connected, and  supercritical regimes, we present 4 regimes where the soft random simplicial complexes $R(n,r,\rho)$ and $C(n,r,\rho)$ exhibit similar behavior to the ones presented by the simplicial complexes $R(n,r)$ and $C(n,r)$ in the mentioned regimes.

\begin{remark}\label{rem6.1}
From (\ref{ineq1}) we have:
\begin{equation*}
\hat{o}_k\leq \beta_k(R(n,r))\leq \hat{o}_k+ f_k^{\geq 2k+3}
\end{equation*}
where $\hat{o}_k$ is the number of components isomorphic to the complex $O_k$ from Definition~\ref{defO_k}, and $f_k^{\geq 2k+3}$ is the number of $k$ faces in components of at least $2k+3$ vertices. Furthermore, by Proposition~\ref{cuentaso_k} the number of $i$-faces of $O_k$ is given by the following expression
$$f_i(O_k)=2^{i+1}\binom{k+1}{i+1}.$$
On the other hand, From (\ref{ineq2}) we have:
\begin{equation*}
\hat{S}_{k+1}\leq \beta_k(C(n,r))\leq \hat{S}_{k+1}+ f_k^{\geq k+3}
\end{equation*}
where $\hat{S}_{k+1}$ is the number of components isomorphic to an empty $k+1$-simplex, and $f_k^{\geq k+3}$ is the number of $k$ faces in components of at least $k+3$ vertices. Furthermore, if $\sigma$ is an empty $k+1$-simplex, the number of $i$-faces of $\sigma$ is given by the following expression
$$f_i(\sigma)=\binom{k+2}{i+1}$$
\end{remark}
\begin{lemma}\label{defr}
Let $R(n,r,\rho)$ be the  soft random  Vietoris-Rips complex and let $C(n,r,\rho)$ be the  soft random \v{C}ech complex with the multiparameter vector $\rho=(p_1,p_2,\dots)$.Then
\begin{enumerate}
\item If $\hat{T}_{k+1}(\rho)$  the number of components isomorphic to an empty $k+1$ simplex and $f_k^{\geq k+3}(\rho)$ is the number of $k$-faces in components of at least $k+3$ vertices, respectively, in $R(n,r,\rho)$. Then 
$$\hat{T}_{k+1}(\rho)\leq \beta_k(R(n,r,\rho))\leq \hat{T}_{k+1}(\rho)+ f_k^{\geq k+3}(\rho)$$

\item If $\hat{S}_{k+1}(\rho)$  the number of components isomorphic to an empty $k+1$ simplex and $f_k^{\geq k+3}(\rho)$ is the number of $k$-faces in components of at least $2k+3$ vertices, respectively, in $R(n,r,\rho)$. Then 
$$\hat{S}_{k+1}(\rho)\leq \beta_k(C(n,r,\rho))\leq \hat{S}_{k+1}(\rho)+ f_k^{k+3\geq}(\rho).$$
\end{enumerate}
\end{lemma}

\begin{proof}
As we discussed in  Remark \ref{rem6.1}, the inequality (\ref{ineq1}) gives us bounds for  $\beta_k(R(n,r))$ in terms of $\hat{o}_k$ and $f_k^{\geq 2k+3}$. The lower bound was obtained from the fact that $O_k$ is the minimal complex required to increase $\beta_k(R(n,r))$ by one. Furthermore, since only components are being considered, each component isomorphic to $O_k$ must increase $\beta_k(R(n,r))$ by one. Thus 
$$\hat{o}_k\leq \beta_k(R(n,r)).$$
However, the possibility of removing $k$-faces in $R(n,r,\rho)$ implies a change in the minimal complex required to increase $\beta_k(R(n,r))$ by one. By definition of the Vietoris-Rips complex  $R(n,r)$, it is not possible to have an empty $k+1$-simplex as a subcomplex of $R(n,r)$ and therefore . Unlike $R(n,r)$, this is not the case for  $R(n,r,\rho)$ where we can achieve an an empty $k+1$-simplex with probability $(1-p_{k+1})$ from a $k+1$-simplex in $R(n,r)$. Considering the above
$$\hat{T}_{k+1}(\rho)\leq \beta_k(R(n,r,\rho)).$$
For the upper bound, notice that any contribution not coming from a component of $R(n,r,\rho)$ isomorphic to an empty $k+1$-simplex, it must  come from a component of at least $k+3$ vertices and at least one $k$-face. Thus
$$\beta_k(R(n,r,\rho))\leq \hat{T}_{k+1}(\rho)+ f_k^{\geq k+3}(\rho).$$
The same argument can be used to prove that
$$\hat{S}_{k+1}(\rho)\leq \beta_k(C(n,r,\rho))\leq \hat{S}_{k+1}(\rho)+ f_k^{k+3\geq}(\rho).$$
\end{proof}

\begin{lemma}\label{xxw1}
Let $s_j$ be the number of subgraphs isomorphic to a certain graph $H_j$ on $k+3$ vertices and $\binom{k+1}{2}+2$ edges, then $f_k^{\geq2k+3}$ can be bound as a finite sum of $s_j$'s.
\end{lemma}
Lemma \ref{xxw1} was taken from the discussion of Theorem 3.1's proof on \cite{MR2770552}.

The following result gives the asymptotic behavior of $E(\beta_k(p))$, for  $C(n,r,\rho)$ and $R(n,r,\rho)$ in the subcritical regime
\begin{theorem} \label{theo6.3}
Let $R(n,r,\rho)$ be the  soft random  Vietoris-Rips complex with the multiparameter vector $\rho=(p_1,p_2,\dots)$. If $\lim\limits_{n\to \infty} \frac{nr^d}{ (1-p_{k+1})\prod\limits_{i=1}^{k}p_i^{\binom{k+1}{i}} }=0$, then 
$$E(\beta_k(R(n,r,\rho)))\sim E(\hat{T}_{k+1}(\rho)).$$
\end{theorem}

\begin{proof}
By Lemma~\ref{defr}:
\begin{equation}\label{ineq6.4}
E(\hat{T}_{k+1}(\rho))\leq E(\beta_k(\rho))\leq E(\hat{T}_{k+1}(\rho))+E( f_k^{\geq k+3}(\rho))
\end{equation}
thus
$$1\leq \frac{E(\beta_k(\rho))}{E(\hat{T}_{k+1}(\rho))}\leq 1+\frac{E( f_k^{\geq k+3}(\rho))}{E(\hat{T}_{k+1}(\rho))}.$$
Now, by Theorem~\ref{newadd}
$$E(\hat{T}_{k+1}(\rho))= (1-p_{k+1})\prod_{i=1}^{k}p_i^{\binom{k+2}{i+1}} E(\hat{T}_{k+1})$$ 
and by  Proposition~\ref{Penrose-2} , $E(\hat{T}_{k+1})\sim n^{k+2}r^{d(k+1)}$. 
Thus by Theorem~\ref{theo6.1}:
$$E(f_k^{k+3\geq}(\rho))\leq c E(s_j(\rho))\leq c_2\prod_{i=1}^{k}p_i^{f_i(s_j(\rho))} E(s_j(\rho)) \leq c_2\prod_{i=1}^{k}p_i^{\binom{k+1}{i+1}}E(s_j(\rho)).$$
Therefore for certain constant $c_3$
$$\frac{E( f_k^{2k+3\geq}(\rho))}{E(\hat{T}_{k+1}(\rho))}\leq\frac{c_3\prod\limits_{i=1}^{k}p_i^{\binom{k+1}{i+1}} E(s_i)}{ (1-p_{k+1})\prod\limits_{i=1}^{k}p_i^{\binom{k+2}{i+1}} E(\hat{T}_{k+1})}$$
and since, by Propositions~\ref{Penrose-1} and \ref{Penrose-2}, $E(s_j)\sim n^{k+3}r^{d(k+2)}$ and $E(\hat{T}_{k+1})\sim n^{k+2}r^{d(k+1)}$ we have
$$\frac{E( f_k^{k+3\geq}(\rho))}{E(\hat{T}_{k+1}(\rho))} \leq\frac{c_3nr^d}{ (1-p_{k+1})\prod\limits_{i=1}^{k}p_i^{\binom{k+1}{i}}}$$
therefore, from our hypothesis we conclude that
$$1\leq \frac{E(\beta_k(\rho))}{E(\hat{T}_{k+1}(\rho))}\leq 1+ \frac{c_3nr^d}{ (1-p_{k+1})\prod\limits_{i=1}^{k}p_i^{\binom{k+1}{i}}}\rightarrow 1$$
and thus 
$$\frac{E(\beta_k(R(n,r\rho)))}{E(\hat{T}_{k+1}(\rho))}\rightarrow 1.\qedhere$$

\end{proof}

The following theorem is the analogue of Theorem~\ref{theo6.3} for the case of the soft random \v{C}ech complex.

\begin{theorem} \label{theo6.9}
Let $C(n,r,\rho)$ the  soft random \v{C}ech complex with respect to a multiparameter vector $\rho=(p_1,p_2,\dots)$. If $\lim\limits_{n\to \infty} \frac{nr^d}{ \prod\limits_{i=1}^{k}p_i^{\binom{k+1}{i}} }=0$, then 
$$E(\beta_k(C(n,r,\rho)))\sim E(\hat{S}_{k+1}(\rho)).$$
\end{theorem}
\begin{proof}
By Lemma~\ref{defr}:
\begin{equation}\label{ineq6.4-}
\hat{S}_{k+1}(\rho)\leq \beta_k(C(n,r,\rho))\leq \hat{S}_{k+1}(\rho)+ f_k^{k+3\geq}(\rho)
\end{equation}
thus
$$E(\hat{S}_{k+1}(\rho))\leq E(\beta_k(C(n,r,\rho)))\leq E(\hat{S}_{k+1}(\rho))+E(f_k^{k+3\geq}(\rho))$$
and 
$$1\leq \frac{E(\beta_k(\rho))}{E(\hat{S}_{k+1}(\rho))}\leq 1+c\frac{E(s_j)}{E(\hat{S}_{k+1}(\rho))}.$$
Now, by Theorem~\ref{theo6.6}
$$E(\hat{S}_{k+1}(\rho))\geq c_1\prod_{i=1}^{k}p_i^{\binom{k+2}{i+1}} n^{k+2}r^{d(k+1)}$$ 
since  $E(\hat{S}_{k+1})\sim n^{k+2}r^{d(k+1)}$.
Thus, by Theorem~\ref{theo6.8}:
$$E(f_k^{k+3\geq}(\rho))\leq  cp_1^2\prod_{i=1}^{k}p_i^{\binom{k+1}{i+1}} n^{k+3}r^{d(k+2)}$$
therefore for certain constant $c_3$
$$\frac{E(f_k^{k+3\geq}(\rho))}{E(\hat{S}_{k+1})}\leq\frac{c_3 p_1^2\prod_{i=1}^{k}p_i^{\binom{k+1}{i+1}} n^{k+3}r^{d(k+2)}}{ \prod_{i=1}^{k}p_i^{\binom{k+2}{i+1}} n^{k+2}r^{d(k+1)}}$$
thus
$$\frac{E(f_k^{k+3\geq}(\rho))}{E(\hat{S}_{k+1})}\leq\frac{c_3  nr^d}{ \prod_{i=1}^{k}p_i^{\binom{k+1}{i}} }$$
therefore, from our hypothesis we conclude that
$$1\leq \frac{E(\beta_k(\rho))}{E(\hat{S}_{k+1}(\rho))}\leq 1+ \frac{c_3  nr^d}{ \prod_{i=1}^{k}p_i^{\binom{k+1}{i}} }\rightarrow 1$$
thus 
$$\frac{E(\beta_k(C(n,r\rho)))}{E(\hat{S}_{k+1}(\rho))}\rightarrow 1.\qedhere$$
\end{proof}

\subsection{Regime of vanishing and non vanishing for $\beta_k$}\label{vantovan}

The following results is the analogue of Theorems 3.9 and 3.10 in \cite{MR2770552}, for the case of the soft random simplicial complexes $R(n,r,\rho)$ and $C(n,r,\rho)$. It describes  the threshold  where the $k$-th homology passes from being trivial to being not trivial.

\begin{theorem} \label{theo6.4}
Let $R(n,r,\rho)$ the  soft random  Vietoris-Rips complex with respect to a multiparameter vector $\rho=(p_1,p_2,\dots)$. Assume that $\lim\limits_{n\to \infty} \frac{nr^d}{ (1-p_{k+1})\prod\limits_{i=1}^{k}p_i^{\binom{k+1}{i}}}=0$.
\begin{enumerate}
\item If   $(1-p_{k+1})\prod\limits_{i=1}^{k}p_i^{\binom{k+2}{i+1}}n^{k+2}r^{d{k+1}}\to 0$, then $\beta_k(R(n,r,\rho))=0$ a.a.s
\item If  $(1-p_{k+1})\prod\limits_{i=1}^{k}p_i^{\binom{k+2}{i+1}}n^{k+2}r^{d{k+1}}\to \infty$, then $\beta_k(R(n,r,\rho))\neq 0$ a.a.s
\end{enumerate}
\end{theorem}
\begin{proof}
\begin{enumerate}
\item Using Markov's inequality Lemma~\ref{Frieze1} and Theorem~\ref{theo6.3}:
$$P(\beta_k(R(n,r,\rho)))\leq E(\beta_k(R(n,r,\rho)))\sim E(\hat{T}_{k+1}(\rho))\rightarrow 0$$
\item  Since
$$E(\hat{T}_{k+1}(\rho))\leq E(\beta_k(\rho))\leq E(\hat{T}_{k+1}(\rho))+E( f_k^{k+3\geq}(\rho))$$
then, by Theorem~\ref{theo6.3}
$$\frac{E(\beta_k(\rho)^2)}{E(\beta_k(\rho))^2}\sim \frac{E(\beta_k(\rho)^2)}{E\hat{T}_{k+1}(\rho))^2}$$
\begin{equation}\label{ineq6.2}
\leq \frac{E(\hat{T}_{k+1}(\rho)^2)}{E(\hat{T}_{k+1}(\rho))^2}+\frac{2E(\hat{T}_{k+1}(\rho))E(f_k^{k+3\geq}(\rho))}{E(\hat{T}_{k+1}(\rho))^2}+\frac{E((f_k^{k+3\geq}(p))^2)}{E(\hat{T}_{k+1}(\rho))^2} 
\end{equation}
let us consider each expression in (\ref{ineq6.2}) individually. For the second expression, notice that
$$\frac{2E(\hat{T}_{k+1}(\rho))E(f_k^{k+3\geq}(\rho))}{E(\hat{T}_{k+1}(\rho))^2}=\frac{2E(f_k^{k+3\geq}(\rho))}{E(\hat{T}_{k+1}(\rho))}\rightarrow 0$$
where the limit is as in the proof of Theorem~\ref{theo6.3}. For the first expression, notice that
$$E(\hat{T}_{k+1}(\rho)^2)=E\left(\left(\sum_{|X|=k+2}\hat{h}_{n,r,k+2,\rho}(X)\right)\left(\sum_{|Y|=k+2}\hat{h}_{n,r,k+2,\rho}(Y)\right)\right)$$
where $\hat{h}_{n,r,k+2,\rho}(X)$ is the characteristic function that is 1 where $R_{n,r,\rho}(X)$ is a component isomorphic to an empty $k+1$-simplex in $R(n,r,\rho)$ and the sum is over all different sets of $k+2$ vertices. Then
\begin{align*}
E(\hat{T}_{k+1}(\rho)^2)&=\sum_{|X|=k+2, |Y|=k+2}E(\hat{h}_{n,r,k+2,\rho}(X)\hat{h}_{n,r,k+2,\rho}(Y)1_{X=Y})\\
 &\qquad+ \sum_{|X|=k+2, |Y|=k+2}E(\hat{h}_{n,r,k+2,\rho}(X)\hat{h}_{n,r,k+2,\rho}(Y)1_{X\neq Y})\\
&=E(\hat{T}_{k+1}(\rho))+ \sum_{|X|=k+2, |Y|=k+2}E(\hat{h}_{n,r,k+2,\rho}(X)\hat{h}_{n,r,k+2,\rho}(Y)1_{X\neq Y})\\
&=E(\hat{T}_{k+1}(\rho))+ \sum_{|X|=k+2, |Y|=k+2}E(\hat{h}_{n,r,k+2,\rho}(X)\hat{h}_{n,r,k+2,\rho}(Y)1_{X\cap Y=\emptyset})\\
&\qquad+ \sum_{|X|=k+2, |Y|=k+2}E(\hat{h}_{n,r,k+2,\rho}(X)\hat{h}_{n,r,k+2,\rho}(Y)1_{|X\cap Y|=l>1})
\end{align*}

now, notice that if $|X\cap Y|=l\geq 1$, then $\hat{h}_{n,r,k+2,\rho}(X)\times\hat{h}_{n,r,k+2,\rho}(Y)=0$, since this means that neither $X$ nor $Y$ span a connected component. Also, $\hat{h}_{n,r,k+2,\rho}(X)$ and $\hat{h}_{n,r,k+2,\rho}(Y)$ are independent whenever $X$ and $Y$ don't share vertices. Thus 
$$E(\hat{T}_{k+1}(\rho)^2)\leq E(\hat{T}_{k+1}(\rho))+ E(\hat{T}_{k+1}(\rho))^2$$
therefore
$$\frac{E(\hat{T}_{k+1}(\rho)^2)}{E(\hat{T}_{k+1}(\rho))^2}\leq \frac{E(\hat{T}_{k+1}(\rho))}{E(\hat{T}_{k+1}(\rho))^2}+ \frac{E(\hat{T}_{k+1}(\rho))^2}{E(\hat{T}_{k+1}(\rho))^2}=\frac{1}{E(\hat{T}_{k+1}(\rho))}+1$$
now, by Theorem~\ref{theo6.2} and the hypothesis, $E(\hat{T}_{k+1}(\rho))$ is going to $\infty$. Then
$$\frac{E(\hat{T}_{k+1}(\rho)^2)}{E(\hat{T}_{k+1}(\rho))^2}\sim 1$$
For the third expression, we may use the same argument. Therefore

$$E((f_k^{k+3\geq}(\rho))^2)=E\left(\left(\sum_{|X|\geq k+3}c_k(X)\hat{h}_{n,r,k+3,\rho}(X)\right)\left(\sum_{|Y|\geq k+3}c_k(Y)\hat{h}_{n,r,k+3,\rho}(Y)\right)\right)$$
where $\hat{h}_{n,r,k+3,\rho}(X)$ is the characteristic function that is $1$ where $R_{n,r,\rho}(X)$ is a component of at least $k+3$ vertices in $R(n,r,\rho)$, $c_k(X)$ is the number of $k$ faces in $R_{n,r,\rho}(X)$. Then
\begin{align*}
E((f_k^{k+3\geq}(\rho))^2)&=\sum_{|X|\geq k+3, |Y|\ \geq k+3}E(c_k(X)c_k(Y)\hat{h}_{n,r,k+3,\rho}(X)\hat{h}_{n,r,k+3,\rho}(Y)1_{X=Y})\\
&\qquad+ \sum_{|X|\geq k+3, |Y|\geq k+3}E(c_k(X)c_k(Y)\hat{h}_{n,r,k+3,\rho}(X)\hat{h}_{n,r,k+3,\rho}(Y)1_{X\neq Y})\\
&=E(f_k^{k+3\geq}(\rho))+ \sum_{|X|\geq k+3, |Y|\geq k+3}E(c_k(X)c_k(Y)\hat{h}_{n,r,k+3,\rho}(X)\hat{h}_{n,r,k+3,\rho}(Y)1_{X\neq Y})\\
&=E(f_k^{k+3\geq}(\rho))+ \sum_{|X|\geq k+3, |Y|\geq k+3}E(c_k(X)c_k(Y)\hat{h}_{n,r,k+3,\rho}(X)\hat{h}_{n,r,k+3,\rho}(Y)1_{X\cap Y=\emptyset})\\
&\qquad+ \sum_{|X|\geq k+3, |Y|\geq k+3}E(c_k(X)c_k(Y)\hat{h}_{n,r,k+3,\rho}(X)\hat{h}_{n,r,k+3,\rho}(Y)1_{|X\cap Y|=l>1})
\end{align*}
now, notice that if $|X\cap Y|=l\geq 1$, then $\hat{h}_{n,r,k+3,\rho}(X)\times\hat{h}_{n,r,k+3,\rho}(Y)=0$, since this means that neither $X$ nor $Y$ span a connected component. Also, $c_k(X)\hat{h}_{n,r,k+3,\rho}(X)$ and $c_k(Y)\hat{h}_{n,r,k+3,\rho}(Y)$ are independent whenever $X$ and $Y$ don't share vertices. Thus 
$$E((f_k^{k+3\geq}(\rho))^2)\leq E(f_k^{k+3\geq}(\rho))+ E(f_k^{k+3\geq}(\rho))^2$$
and thus
$$\frac{E((f_k^{k+3\geq}(\rho))^2)}{E(\hat{T}_{k+1}(\rho)))^2}\leq \frac{E(f_k^{k+3\geq}(\rho))}{E(\hat{T}_{k+1}(\rho))^2}+ \frac{E(f_k^{k+3\geq}(p))^2}{E(\hat{T}_{k+1}(\rho))^2}.$$
Since, by the hypothesis, $E((f_k^{k+3\geq}(\rho)) /E(\hat{T}_{k+1}(\rho))\rightarrow 0$ and $1/E(\hat{T}_{k+1}(\rho))\rightarrow 0$, it follows that 
$$\frac{E((f_k^{k+3\geq}(\rho))^2)}{E(\hat{T}_{k+1}(\rho))^2}\rightarrow 0. $$
Combining the above with Equation (\ref{ineq6.2}), we have that
$$\lim_{n\to \infty}\frac{E(\beta_k(\rho)^2)}{E(\beta_k(\rho))^2}\leq 1.$$
On the other hand, since  $0\leq var(\beta_k(\rho))= E(\beta_k(\rho)^2)- E(\beta_k(\rho))^2$, then 
$$1\leq \lim_{n\to\infty} \frac{E(\beta_k(\rho)^2)}{E(\beta_k(\rho))^2}\leq 1.$$
Finally, by the second moment method Lemma~\ref{Frieze4}:
$$P(\beta_k(\rho)=0)\leq \frac{E(\beta_k(\rho)^2)}{E(\beta_k(\rho))^2}-1\rightarrow 0$$
\end{enumerate}
\end{proof}

The following theorem is the analogue of Theorem~\ref{theo6.4} for the case of the soft random \v{C}ech complex.

\begin{theorem} \label{theo6.10}
Let $C(n,r,\rho)$ the  soft random \v{C}ech complex with respect to a multiparameter vector $\rho=(p_1,p_2,\dots)$. Assume that $\lim\limits_{n\to \infty} \frac{nr^d}{ \prod_{i=1}^{k}p_i^{\binom{k+1}{i}} }=0$.
\begin{enumerate}
\item If   $\prod\limits_{i=1}^{k}p_i^{\binom{k+2}{i+1}}n^{k+2}r^{d{k+1}}\to 0$, then $\beta_k(C(n,r,\rho))=0$ a.a.s
\item If  $\prod\limits_{i=1}^{k}p_i^{\binom{k+2}{i+1}}n^{k+2}r^{d{k+1}}\to \infty$, then $\beta_k(C(n,r,\rho))\neq 0$ a.a.s
\end{enumerate}
\end{theorem}
\begin{proof}
\begin{enumerate}
\item[]
\item Using Markov's inequality Lemma~\ref{Frieze1} and Theorem~\ref{theo6.3}:
$$P(\beta_k(C(n,r,\rho)))\leq E(\beta_k(C(n,r,\rho)))\sim E(\hat{S}_{k+1}(\rho))\rightarrow 0$$
\item  Since
$$\hat{S}_{k+1}(\rho)\leq \beta_k(C(n,r,\rho))\leq \hat{S}_{k+1}(\rho)+ f_k^{k+3\geq}(\rho)$$
then, by Theorem~\ref{theo6.9}
$$\frac{E(\beta_k(\rho)^2)}{E(\beta_k(\rho))^2}\sim \frac{E(\beta_k(\rho)^2)}{E(\hat{S}_{k+1}(\rho))^2}$$
\begin{equation}\label{ineq6.2-}
\leq \frac{E(\hat{S}_{k+1}(\rho)^2)}{E(\hat{S}_{k+1}(\rho))^2}+\frac{2E(\hat{S}_{k+1}(\rho))E(f_k^{k+3\geq}(\rho))}{E(\hat{S}_{k+1}(\rho))^2}+\frac{E((f_k^{k+3\geq}(\rho))^2)}{E(\hat{S}_{k+1}(\rho))^2} 
\end{equation}
let us consider each expression in (\ref{ineq6.2-}) individually. For the second expression, notice that

$$\frac{2E(\hat{S}_{k+1}(\rho))E(f_k^{k+3\geq}(\rho))}{E(\hat{S}_{k+1}(\rho))^2}=\frac{2E(f_k^{k+3\geq}(\rho))}{E(\hat{S}_{k+1}(\rho))}\rightarrow 0$$
where the limit is as in the proof  of Theorem~\ref{theo6.9}. For the First expression, notice that
$$E(\hat{S}_{k+1}(\rho)^2)=E\left(\left(\sum_{|X|=k+2}\hat{h}_{n,r,k+2,\rho}(X)\right)\left(\sum_{|Y|=k+2}\hat{h}_{n,r,k+2,\rho}(Y)\right)\right)$$
where $\hat{h}_{n,r,k+2,\rho}(X)$ is the characteristic function that is 1 where $C_{n,r,\rho}(X)$ is a component isomorphic to an empty $k+1$-simplex in $C(n,r,\rho)$ and the sum is over all different sets of $k+2$ elements. Then, similar to the calculations made for $R(n,r,\rho)$:
$$E(\hat{S}_{k+1}(\rho)^2)\leq E(\hat{S}_{k+1}(\rho))+ E(\hat{S}_{k+1}(\rho))^2$$
therefore
$$\frac{E(\hat{S}_{k+1}(\rho)^2)}{E(\hat{S}_{k+1}(\rho))^2}\leq \frac{E(\hat{S}_{k+1}(\rho))}{E(\hat{S}_{k+1}(\rho))^2}+ \frac{E(\hat{S}_{k+1}(\rho))^2}{E(\hat{S}_{k+1}(\rho))^2}=\frac{1}{E(\hat{S}_{k+1}(\rho))}+1$$
now, by Theorem~\ref{theo6.2} and the hypothesis, $E(\hat{S}_{k+1}(\rho))$ is going to $\infty$. Then
$$\frac{E(\hat{S}_{k+1}(\rho)^2)}{E(\hat{S}_{k+1}(\rho))^2}\sim 1$$
For the third expression, we may use the same argument used  in the proof of Theorem~\ref{theo6.4}. Therefore
$$E((f_k^{k+3\geq}(\rho))^2)\leq E(f_k^{k+3\geq}(\rho))+ E(f_k^{k+3\geq}(\rho))^2$$
and thus
$$\frac{E((f_k^{k+3\geq}(\rho))^2)}{E(\hat{S}_{k+1}(\rho))^2}\leq \frac{E(f_k^{k+3\geq}(\rho))}{E(\hat{S}_{k+1}(\rho))^2}+ \frac{E(f_k^{k+3\geq}(p))^2}{E(\hat{S}_{k+1}(\rho))^2}$$
since under the hypothesis, $E((f_k^{k+3\geq}(\rho)) /E(\hat{S}_{k+1}(\rho)\rightarrow 0$ and $1/E(\hat{S}_{k+1}(\rho)\rightarrow 0$, it follows that
$$\frac{E((f_k^{k+3\geq}(\rho))^2)}{E(\hat{S}_{k+1}(\rho))^2}\rightarrow 0 $$
thus
$$\lim_{n\to \infty}\frac{E(\beta_k(\rho)^2)}{E(\beta_k(\rho))^2}\leq 1.$$
On the other hand, since  $0\leq var(\beta_k(\rho))= E(\beta_k(\rho)^2)- E(\beta_k(\rho))^2$, then
$$1\leq \lim_{n\to \infty}\frac{E(\beta_k(\rho)^2)}{E(\beta_k(\rho))^2}\leq 1.$$
Finally, by the second moment method Lemma~\ref{Frieze4}:
$$P(\beta_k(\rho)=0)\leq \frac{E(\beta_k(\rho)^2)}{E(\beta_k(\rho))^2}-1\rightarrow 0$$
\end{enumerate}
\end{proof}

\section{The critical regime}\label{section4}

We now analyze the critical regime $nr^d\rightarrow \rho\in (0,\infty)$, where the expected value of the Betti numbers increase linearly with the number of points $n$.
The following result is the analogue of Theorem 4.1 in \cite{MR2770552}, for the case of the soft random simplicial complexes $R(n,r,\rho)$ and $C(n,r,\rho)$.
The proofs follow the strategy used in \cite{MR2770552} for the corresponding results on random geometric Vietoris-Rips and random geometric \v{C}ech complexes.
\begin{theorem}\label{theo4.1} Consider the soft random simplicial complexes $R(n,r,\rho)$ and $C(n,r,\rho)$. Let $\epsilon>0$ be given. Assume that $nr^d \rightarrow \lambda\in (0,\infty)$. Then
\begin{enumerate}
\item  If $(1-p_{k+1})\prod\limits_{i=1}^{k}p_i^{\binom{k+2}{i+1}}\geq \frac{1}{n^{\epsilon}}$, 
then $E(\beta_k(R(n,r,\rho)))=O(n)$ and $E(\beta_k(R(n,r,\rho)))=\Omega(n^{1-\epsilon})$
\item If $\prod\limits_{i=1}^{k}p_i^{\binom{k+2}{i+1}}\geq \frac{1}{n^{\epsilon}}$, 
then $E(\beta_k(C(n,r,\rho)))=O(n)$ and $E(\beta_k(C(n,r,\rho)))=\Omega(n^{1-\epsilon})$
\end{enumerate}
\end{theorem}
\begin{proof} From (\ref{ineq6.4}) it follows that
$$E(\hat{T}_{k+1}(\rho))\leq E(\beta_k(R(n,r,\rho)))\leq E(\hat{T}_{k+1}(\rho))+ E(f_k^{k+3\geq}(\rho))$$
by Theorems~\ref{newadd} and \ref{theo6.2}, for a constant $c_1$:
$$(1-p_{k+1})\prod_{i=1}^{k}p_i^{\binom{k+2}{i+1}}E(\hat{T}_{k+1})\leq E(\beta_k(R(n,r,\rho)))\leq (1-p_{k+1}) \prod_{i=1}^{k}p_i^{\binom{k+2}{i+1}} E(\hat{T}_{k+1})+ c_1 \prod_{i=1}^{k}p_i^{\binom{k+1}{i+1}} E(s_j)$$
therefore for a constant $c$
$$(1-p_{k+1})\prod_{i=1}^{k}p_i^{\binom{k+2}{i+1}} E(\hat{T}_{k+1})\leq E(\beta_k(R(n,r,\rho)))\leq c (E(\hat{T}_{k+1})+  E(s_j))$$
then, by the hypothesis and Proposition~\ref{Penrose-3} and \ref{Kahle-1} there exist constants $d_1$ and $d_2$ such that $E(\hat{T}_{k+1})\sim d_1n $ and $E(s_j)\sim d_2n$. So there are constants $c'_1$ and $c'$ such that
$$c_1' n^{1-\epsilon}\leq E(\beta_k(R(n,r,\rho)))\leq c' n.$$
The proof is analogous for $E(\beta_k(C(n,r,\rho)))$.
\end{proof}

\section{The supercritical regime}\label{section5}

We now study the supercritical regime $nr^d\rightarrow \infty$. In the super critical regime, following Kahle\cite{MR2770552}, we consider $n$ points $X_1,\dots,X_n$ uniformly distributed in a  bounded convex body $K$, since as noted in \cite{MR2770552}, some assumption on density must be made in order to deduce topological results in  denser regimes. We recall  that a  bounded convex body  is defined to be a compact, convex set with nonempty interior. The following theorem gives  the behavior of the $\beta_k(\rho)$ in the supercritical regime.

The main tool used to prove the following theorem will be discrete Morse theory. Discrete Morse theory is a powerful tool used to analyze the topology and geometry  of cellular complexes and pseudomanifolds. It provides information about the topological structure  of these spaces from functions defined on them.

We now introduce some terminology of discrete Morse
theory and state the main theorem which we will use. For information about the subject we refer the reader to \cite{MR1939695}.

\begin{definition}
\begin{enumerate}\label{defA03}
\item[]
\item Let $\sigma$ and $\beta$ two faces a of a simplicial complex. We write $\{\sigma_1 \prec\beta_1\}$ iff $\sigma=\text{dimension }\beta - 1$
\item A discrete vector field $V$ of a simplicial complex $X$ is a collection of pairs of faces $\{\sigma_1 \prec\beta_1\}$  such that each face is in at most  one pair.
\item Given a discrete vector field $V$ for simplicial complex $X$, a closed path in $V$ is a sequence of faces 
$$\sigma_1\prec\beta_1\succ\sigma_2\prec\beta_2\succ\sigma_3\prec\beta_3\succ\dots\sigma_n\prec\beta_n\succ\sigma_{n+1}$$
with $\sigma_i\neq \sigma_{i+1}$, every pair $\{\sigma_i\prec\beta_i\}\in V$  and $\sigma_1=\sigma_{n+1}$. Notice that, $\{\beta_i>\sigma_{i+1}\}\not\in V$ since each face is in at most one pair. 
\item A vector field $V$ over a simplicial complex $X$ is  called gradient vector field if there are no closed paths in $V$.
\item A  $k$-face $\sigma \in X$ is called critical with respect to the vector field $V$ iff it is not matched in $V$.
\end{enumerate}
\end{definition}

\begin{theorem}[Fundamental theorem of discrete Morse theory]\label{morse}
Suppose the pair $(X,V)$ is a simplicial complex together with a discrete gradient vector field $V$. Then $X$ is homotopy equivalent to a CW complex with one cell of dimension $k$ for each  $k$-dimensional face $\sigma$ of $X$ which is critical with respect to $V$.
\end{theorem}

The following result is the analogue of Theorem 5.1 in \cite{MR2770552}, for the case of the soft random Vietoris-Rips complex $R(n,r,\rho)$.

\begin{theorem}\label{theo5.1}
Let $R(n,r,\rho)$ be the soft random  Vietoris-Rips complex for which the  $n$ points $ \{X_1,\dots,X_n\}$ are taken i.i.d uniformly from a  bounded convex body $K \subset \mathbb{R}^d$. Suppose $r=\omega(n^{-\frac{1}{d}})$ and write $W=nr^d$. Then
$$E( \beta_k(R(n,r,\rho)))=O\left(nW^k\prod_{i=1}^{k}p_i^{\binom{k+1}{i+1}}e^{-cW\prod\limits_{i=1}^{k+1}p^{\binom{k+1}{i}}}\right).$$
\end{theorem}

\begin{proof}
With probability 1, no two points are at the same distance from the origin. We therefore assume that
$$||X_1||<||X_2||<\dots<||X_n||.$$
Let $V$ the vector field for the soft random Vietoris-Rips complex $R(n,r,\rho)$ which pair every  $k$-face $F=(X_{i_1}, \dots, X_{i_{k+1}})$ with a face $\{X_{i_0}\}\cup F$ for $i_0<i_1$ and as small as possible. It is proved in Theorem 5.1 in \cite{MR2770552} that $V$ is a gradient vector field (see Definition~\ref{defA03}). Thus, we can use Theorem \ref{morse} the fundamental theorem of the discrete Morse theory to guarantee that $E(\beta_k(R(n,r,\rho)))\leq E(C_k)$, where $C_k$ is the number of critical $k$-faces under the vector field $V$.

The probability $Q_1$ that a set of $k+1$ vertices span a $k$-face in $R(n,r,\rho)$ satisfies that 
$$Q_1=O(\prod\limits_{i=1}^{k}p^{\binom{k+1}{i+1}}r^{dk})$$
since if $A$ is a set of $k+1$ elements: 
$$P(\text{$A$ spans a k face in $R(n,r,\rho)$}|\text{$A$ spans a $k$-face in $R(n,r)$})=\prod_{i=1}^{k}p^{\binom{k+1}{i+1}}$$
and by Proposition~\ref{Penrose-1}
$$P(  \text{$A$ spans a $k$-face in $R(n,r)$}  )=O(r^{dk}).$$
Now, notice that, by definition of the pairing relation in $V$, a $k$-face $F=(X_{i_1}, \dots, X_{i_{k+1}})$ with $j<l \implies i_j < i_l$ is paired in the vector field $V$ iff
\begin{enumerate}
\item If there is  a vertex $X_a$ that is neighbor in $R(n,r,\rho)$  of each vertex in $F$  and  such that $a<i_1$ then $\{F\prec F\cup\{X_a\}\}\in V$.
\item If there is no  vertex $X_b$ that is a neighbor in $R(n,r,\rho)$  of each vertex in $F-\{x_{i_1}\}$ and  such that $b<i_1$ then $\{(X_{i_2}...X_{i_{k+1}})\prec (X_{i_2}...X_{i_{k+1}})\cup\{X_1\}=F\}\in V$.
\end{enumerate}
Thus, if we assume that $F=(X_{i_1}, \dots, X_{i_{k+1}})$ is a critical face with $j<l \implies i_j < i_l$, then there exists a point $X_{i_0}$ with $i_0<i_1$ that is a neighbor of $X_{i_2},...,X_{i_{k+1}}$ in $G(n,r,\rho)$. Furthermore, notice that $X_{i_0}$ and $F$ can't span a $k+1$-face since in this case $F$ would be paired.

By Lemma~\ref{interccc}, for the set
$$I=\bigcap_{j=1}^{k+1}B(X_{i_j},r)\cap B(0, ||X_{i_1}||)$$
we have that $Leb(I)\geq \epsilon_d r^d$ for some constant $\epsilon_d >0$. Suppose that a vertex $y\notin\{X_{i_0},\dots,X_{i_{k+1}}\}$ falls in the region $I$. Then, $F$ would be paired  in the discrete vector field  $V$ to a $(k+1)$-face $\tilde{F}_y \coloneqq F \cup \{y\}$ with probability  $\prod\limits_{i=1}^{k+1}p_i^{\binom{k+1}{i}}$. Therefore, if $y$  runs over all of the other $n-k-2$ vertices and $Q_2=P(\text{$F$ is critical}\mid \text{$F$ is a $k$-face})$ we have 

\begin{align*}
 Q_2 &\leq P(\forall y \in \mathcal{X}_n\backslash\{X_{i_0},\dots, X_{i_{k+1}}\}, y\not\in I) \\
 & + P(\forall y \in \mathcal{X}_n\backslash\{X_{i_0},\dots, X_{i_{k+1}}\}, y\in I, \text{$\tilde{F}_y$ is not a $(k+1)$-face in $R(n,r,\rho)$})
\end{align*}
now, since the points are uniformly distributed over $K$:
$$P(\text{$F$ is critical}\mid \text{$F$ is a $k$-face})\leq \left(1-\frac{Leb(I)}{Leb(K)}\right)^{n-k-2}+ \left(\frac{Leb(I)}{Leb(K)}(1-\prod\limits_{i=1}^{k+1}p_i^{\binom{k+1}{i}})\right)^{n-k-2}$$

$$\leq e^{(n-k-2)cr^d}+ e^{-(n-k-2)r^d\prod\limits_{i=1}^{k+1}p_i^{\binom{k+1}{i}}}$$
$$\leq 2 e^{-c(n-k-2)r^d\prod\limits_{i=1}^{k+1}p_i^{\binom{k+1}{i}}}$$
By the fundamental theorem of the discrete Morse theory, Theorem \ref{morse}, we can bound the $k$-th Betti number by the number of  critical $k$-faces of the vector field $V$. Thus if $C_k$ is the number of  critical $k$-faces of the vector field $V$:
$$E(\beta_k(R(n,r,\rho)))\leq E(C_k)=\binom{n}{k+1}Q_1Q_2=O\left( \prod\limits_{i=1}^{k}p_i^{\binom{k+1}{i+1}}n^{k+1}r^{dk}e^{-c(n-k-2)r^d\prod\limits_{i=1}^{k+1}p_i^{\binom{k+1}{i}}}\right).$$
Thus, since $r^d\to 0$  and $W=nr^d$ by hypothesis, we have
$$E(\beta_k(R(n,r,\rho)))=O\left(\prod\limits_{i=1}^{k}p_i^{\binom{k+1}{i+1}}n W^k e^{-cW\prod\limits_{i=1}^{k+1}p_i^{\binom{k+1}{i}}}\right)$$
\end{proof}

\begin{corollary}
Let $R(n,r,\rho)$ be the soft Vietoris-rips complex for which the  $n$ points were taken i.i.d uniformly from a  bounded convex body $K$ in $\mathbb{R}^d$. Write $W=nr^d$ and suppose $W\prod\limits_{i=1}^{k+1}p_i^{\binom{k+1}{i}}\rightarrow \infty$. Then
$$E( \beta_k(R(n,r,\rho)))=o(n)$$
\end{corollary}
\begin{proof}
By  Theorem~\ref{theo5.1}, for a constant $\lambda$:
$$E(\beta_kR(n,r,\rho))\leq  \lambda \prod\limits_{i=1}^{k}p^{\binom{k+1}{i+1}}n W^k e^{-cW\prod\limits_{i=1}^{k+1}p_i^{\binom{k+1}{i}}}\leq  \lambda n W^k e^{-cW\prod\limits_{i=1}^{k+1}p_i^{\binom{k+1}{i}}}$$
thus
$$\frac{E(\beta_kR(n,r,\rho))}{n}\leq  \lambda \frac{ W^k}{ e^{cW\prod\limits_{i=1}^{k+1}p_i^{\binom{k+1}{i}}}}\rightarrow 0.$$
\end{proof}

\section{The connected  regime}\label{section6}
\begin{definition}
A topological space $X$ is $k$-connected if every pointed map from the $i$-dimensional sphere $S^i\rightarrow X$ is homotopic to the constant map i.e. $\pi_i(X)$ is trivial for all $0\leq i\leq k$.
\end{definition}

The following result is the analogue of Theorem 6.1 in \cite{MR2770552}, for the case of the soft random Vietoris-Rips complex $R(n,r,\rho)$.
\begin{theorem}\label{theo666}
Let $R(n,r,\rho)$ be the soft Vietoris-rips for which the  $n$ points were taken i.i.d uniformly from a  bounded convex body $K$ in $\mathbb{R}^d$. Fix $k\geq 0$, let $r\geq c_k(\frac{\log(n)}{n})^{\frac{1}{d}}$ , and suppose that for each $i\leq k$, $p_i(n)$ is a sequence bounded away from zero. Then $R(n,r,\rho)$ is $k$-connected a.a.s.
\end{theorem}

\begin{proof} 
Notice that, since the property of $k$-connectivity is monotone respect to $r$, it is enough to prove that the theorem is true for $r= c_k(\frac{\log(n)}{n})^{\frac{1}{d}}$, thus we assume that  $r= c_k(\frac{\log(n)}{n})^{\frac{1}{d}}$. By Theorem~\ref{theo5.1}:
$$E(\beta_k(R(n,r,\rho)))=O(\prod\limits_{i=1}^{k}p^{\binom{k+1}{i+1}}n W^k e^{-cW\prod\limits_{i=1}^{k+1}p^{\binom{k+1}{i}}}),$$
where we recall that $W=nr^d$. Therefore, there exists a constant $\lambda$ such that
\begin{align*}
E(\beta_k(R(n,r,\rho)))&\leq \lambda \prod\limits_{i=1}^{k}p^{\binom{k+1}{i+1}}n W^k e^{-cW\prod\limits_{i=1}^{k+1}p^{\binom{k+1}{i}}}\\
&\leq\lambda n (nr^d)^k e^{-cnr^d\prod\limits_{i=1}^{k+1}p^{\binom{k+1}{i}}}\\
&\leq \lambda n (c_k^d\log(n))^k e^{-cc_k^d\log(n)\prod\limits_{i=1}^{k+1}p^{\binom{k+1}{i}}}\\
&\leq  \lambda (c_k^d\log(n))^k n^{1-cc_k^d\prod\limits_{i=1}^{k+1}p^{\binom{k+1}{i}}}\\
&\leq  \lambda\frac{ (c_k^d\log(n))^k}{ n^{cc_k^d\prod\limits_{i=1}^{k+1}p^{\binom{k+1}{i}}-1}}
\end{align*}
by choosing $c_k^d>\frac{1}{c\prod\limits_{i=1}^{k+1}p^{\binom{k+1}{i}}}$, we have that 
$$E(\beta_k)\rightarrow 0.$$ 
Also, for $l<k$, by Theorem~\ref{theo5.1}:
$$E(\beta_lR(n,r,\rho))=O( \prod\limits_{i=1}^{l}p^{\binom{l+1}{i+1}}n W^l e^{-cWp^{l+1}})$$
$$\leq  \lambda\frac{ (c_k^d\log(n))^l}{ n^{cc_k^d\prod\limits_{i=1}^{l+1}p^{\binom{l+1}{i}}-1}}\rightarrow 0$$
since  $c_k^d>\frac{1}{c\prod\limits_{i=1}^{k+1}p^{\binom{k+1}{i}}}>\frac{1}{c\prod\limits_{i=1}^{l+1}p^{\binom{l+1}{i}}}$.
\end{proof}
The following result  ensures that the random \v{C}ech complex $C(n,r)$ is contractible in the connected regime. 
\begin{theorem}[Theorem 6.1 \cite{MR2770552} ]\label{theo5.5}
Let $C(n,r)$ be the random \v{C}ech complex for which the  $n$ points were taken i.i.d uniformly from a  bounded convex body $K$ in $\mathbb{R}^d$. Then, there exists a constant $c$ depending on $K$, such that if $r\geq c(\frac{\log(n)}{n})^{\frac{1}{d}}$, then $C(n,r)$ is contractible a.a.s.
\end{theorem}

Theorem~\ref{theo5.5} was proved in \cite{MR2770552} using the Nerve Theorem and showing that in the connected regime, the balls $B(X_i, \frac{r}{2})$ cover $K$. In this case, however, the possibility of removing faces in $C(n,r,\rho)$ does not  allow us to use this argument since the simplicial complex $C(n,r,\rho)$ is no longer the nerve of a cover. When the probabilities in the probability vector $\rho$ are sufficiently high, however, we may show that $C(n,r,\rho)$ shows the same asymptotic behavior as $C(n,r)$, which we state in the following theorem.

\begin{theorem}\label{theo5.6}
Let $\epsilon>0$ be given and let $C(n,r,\rho)$ be the soft random \v{C}ech complex for which  $n$ points were taken i.i.d uniformly from a  bounded convex body $K$ in $\mathbb{R}^d$.  If $p_i\geq 1-\frac{1}{n^{2+i+\epsilon}}$ for all $i$, then there exists a constant $c_K$ depending on $K$, such that if $r\geq c_K(\frac{\log(n)}{n})^{\frac{1}{d}}$, then $C(n,r,\rho)$ is contractible a.a.s.
\end{theorem}

\begin{proof}
By Theorem~\ref{theo5.5}, $C(n,r)$ is contractible a.a.s. Therefore, with probability tending to 1, $\beta_k(C(n,r))=0$ for every $k$. On the other hand, by definition,  if $C(n,r,\rho)$ keeps every possible face of $C(n,r)$, then $C(n,r,\rho)=C(n,r)$. Thus:
\begin{align*}
\prod_{i=1}^{n}p_i^{f_i(C(n,r))}=P(C(n,r,\rho)\text{ keeps every face of }C(n,r))&= P(C(n,r,\rho)=C(n,r))\\
&\leq P(\beta_k(C(n,r,\rho))=\beta_k(C(n,r))).
\end{align*}
On the other hand
$$\prod_{i=1}^{n}p_i^{f_i(C(n,r))}\geq \prod_{i=1}^{n}p_i^{\binom{n}{i+1}}\geq \prod_{i=1}^n(1-\frac{1}{n^{2+i+\epsilon}})^{\binom{n}{i+1}}\sim \prod_{i=1}^ne^{-\frac{c}{n^{1+\epsilon}}}\sim e^{-\frac{cn}{n^{1+\epsilon}}}\to 1$$
therefore $p(\beta_k(C(n,r,\rho))=\beta_k(C(n,r)))\rightarrow 1$ and $\beta_k(C(n,r,\rho))=0$ a.a.s.
\end{proof}

\subsection{Non vanishing to vanishing threshold}
In Section~\ref{vantovan}, we gave a threshold for the appearance of non-trivial Betti numbers. We present here a second threshold  where the $k$-dimensional homologies of the random simplicial complexes $R(n,r,\rho)$ and $C(n,r,\rho)$ vanish again. 
\begin{lemma}[Lemma 6.6 \cite{MR2770552}]\label{lemma5.7}
Suppose $H$ is a feasible subgraph, $r=\Omega(n^{-\frac{1}{d}})$ and $r=o(\frac{(\log(n)}{n}^{\frac{1}{d}}))$. The the geometric random graph $G(n,r)$ has at least one connected component isomorphic to $H$ a.a.s.
\end{lemma}
The previous Lemma is stated in \cite{MR2770552} without proof. However, as they clarify, this result can be proved using the techniques in Chapter 3 of  Penrose's  book  \cite{MR1986198}, the second moment method and  the fact that in this regime  $E(J_{n,r,A}(\Gamma))$ dominates $var(J_{n,r,A}(\Gamma))$. In the following lemma, we use the above to show the analogue of Lemma~\ref{lemma5.7} for the case of the soft random graph $G(n,r,p)$.

\begin{lemma}\label{lemmma5.8}
Suppose $H$ is a feasible subgraph, $r=\Omega(n^{-\frac{1}{d}})$ and $r=o(\frac{(\log(n)}{n}^{\frac{1}{d}}))$. Then, if $p\rightarrow 1$, the soft random graph $G(n,r,p)$ has at least one connected component isomorphic to $H$ a.a.s.
\end{lemma}
\begin{proof}
Let $H$ be a  fixed feasible subgraph. Since we have with probability tending to 1 that there exists a component isomorphic to $H$ in $G(n,r)$,  this component survives in $G(n,r,p)$ with probability $p^{f_1(H)}$. Thus, if $J_{n,r}(H) $ and $J_{n,r}(H)(p)$ count the number of components isomorphic to $H$ in $G(n,r)$ and $G(n,r,p)$, respectively
$$P(J(H)(p)>0 \mid J(H)>0)\geq p^{f_1(H)}.$$
By Lemma~\ref{lemma5.7}, $P(J(H)>0)\rightarrow 1$, and therefore
\begin{align*}
P(J(H)(p)>0)\geq P(J(H)(p)>0,J(H)>0)&=P(J(H)(p)>0|J(H)>0)P(J(H)>0)\\
&\geq p^{f_1(H)}P(J(H)>0)\rightarrow 1.
\end{align*}

\end{proof}

\begin{theorem}\label{need1}
Let $R(n,r,\rho)$ be the soft Vietoris-rips complex for which  the  $n$ points were taken i.i.d uniformly from a  bounded convex body $K \subset \mathbb{R}^d$. For any fixed $k\geq 1$:
\begin{enumerate}
\item If either $r^{d{k+1}}=o((1-p_{k+1})^{-1}\prod\limits_{i=1}^{k}p_i^{-\binom{k+2}{i+1}}n^{-(k+2)})$ or else $r= \omega((\frac{\log(n)}{n})^{\frac{1}{d}})$ and  for all $i\leq k$, $p_i(n)$ is a sequence bounded away from zero, then $H_k=0$ a.a.s.
\item If $r^{d{k+1}}=\omega((1-p_{k+1})^{-1}\prod\limits_{i=1}^{k}p_i^{-\binom{k+2}{i+1}}n^{-(k+2)})$, $r= o((\frac{\log(n)}{n})^{\frac{1}{d}})$, and for all $i\leq k$, $p_i(n)\to 1$, then $H_k\neq 0$ a.a.s.
\end{enumerate}
\end{theorem}
\begin{proof}
Part 1 comes from the results in the sub-critical and  connected regime, respectively. For part 2, notice that, from Lemma~\ref{defr}, if $\hat{T}_{k+1}(\rho)$ is the number of components isomorphic to an empty $k+1$-simplex in $R(n,r,\rho)$, then 
$$\hat{T}_{k+1}(\rho)\leq \beta_k(R(n,r,\rho)).$$
On the other hand, by Lemma~\ref{lemmma5.8}, with probability tending to 1 
$R(n,r,\rho)$ has a component isomorphic to an empty $k+1$-simplex. Thus, a.a.s
$$1\leq \hat{T}_{k+1}(\rho)\leq \beta_k(R(n,r,\rho)).$$
\end{proof}
\begin{theorem}\label{need2}
Let $C(n,r,\rho)$ be the soft \v{C}ech complex for which the  $n$ points were taken i.i.d uniformly from a bounded convex body $K \in \mathbb{R}^d$. For any fixed $k\geq 1$
\begin{enumerate}
\item If either $r^{d{k+1}}=o(\prod\limits_{i=1}^{k}p_i^{-\binom{k+2}{i+1}}n^{-(k+2)})$ or else $r= \omega((\frac{\log(n)}{n})^{\frac{1}{d}})$ and for all $i$, $p_i\geq 1-\frac{1}{n^{2+i+\epsilon}}$, then $H_k=0$ a.a.s.
\item If $r^{d{k+1}}=\omega(\prod\limits_{i=1}^{k}p_i^{-\binom{k+2}{i+1}}n^{-(k+2)})$, $r= o((\frac{\log(n)}{n})^{\frac{1}{d}})$, and for all $i\leq k$, $p_i(n)\to 1$, then $H_k\neq 0$  a.a.s.
\end{enumerate}
\end{theorem}
\begin{proof}
Part 1 comes from the results in the sub-critical and  connected regime, respectively. As in the case of $R(n,r,\rho)$, Part 2 follows from Lemma~\ref{defr} and Lemma~\ref{lemmma5.8}. The argument being that asymptotically, we can get a component that is isomorphic to an empty 
$k+1$-simplex. Thus 
$$1\leq \hat{S}_{k+1}(\rho)\leq \beta_k(C(n,r,\rho)).$$

\end{proof}

\section*{Acknowledgements}
I want to thank my advisor Antonio Rieser for guiding me throughout my research and for sharing his knowledge with me. Additionally, I would like to thank  Octavio Arizmendi and Arturo Jaramillo for helpful comments on earlier versions of the manuscript.

\appendix

\section{General Probability results}

\begin{lemma}[Lemma 20.1 \cite{MR3675279}]\label{Frieze1} Let $X$ be  a non-negative random variable. Then, for all $t>0$:
$$P(X\geq t)\leq \frac{E(X)}{t}.$$
\end{lemma}

\begin{lemma}[Lemma 20.2 \cite{MR3675279}]\label{Frieze2} Let $X$ be  a non-negative integer valued random variable. Then:
$$P(X>0)\leq E(X).$$
\end{lemma}

\begin{lemma}[Lemma 20.3 \cite{MR3675279}]\label{Frieze3} If $X$ is a random variable with a finite mean and variance. Then, for all $t>0$:
$$P(|X-E(X)|\geq t)\leq \frac{var(X)}{t^2}.$$
\end{lemma}

\begin{lemma}[Lemma 20.4 \cite{MR3675279}]\label{Frieze4} If X is a non-negative integer valued random variable then
$$P(X=0)\leq \frac{var(X)}{(E(X))^2}=\frac{E(X^2)}{(E(X))^2}-1.$$
\end{lemma}

\begin{lemma}[Lemma 5.3  \cite{MR2770552}]\label{interccc} There is a constant $\epsilon_d>0$ such that the following holds for every $r>0$. Let $l\geq 1$ and $\{y_0,\dots,y_l\}\subset \mathbb{R}^d$ be an $(l+1)$-tuple of points, such that
$$||y_0||\leq ||y_1||\leq \dots \leq ||y_l||$$
and $\frac{1}{2}r\leq ||y_1||$. If $||y_0-y_1||>r$ and $||y_i-y_j||\leq r$ for every other $0\leq i<j\leq l$, then the intersection
$$I=\bigcap_{i=1}^l B(y_i,r)\cap B(0,||y_1||) $$
satisfies $Leb(I)\geq \epsilon_dr^d$.
\end{lemma}

\section*{Declarations}

\subsection{Ethical Approval }
Not applicable
\subsection{Funding }
This work was supported by the the CONAHCYT
postgraduate studies scholarship number 940355 and by the grant
N62909-19-1-2134 from the US Office of Naval Research Global and the Southern Office of Aerospace Research and Development of the US Air Force of Scientific Research.

\subsection{Availability of data and materials }
Not applicable

\end{document}